\documentclass[11pt, reqno]{amsart}
\usepackage{latexsym}
\usepackage{amssymb}
\usepackage{amsmath}
\usepackage{amsfonts}
\usepackage{mathrsfs}
\usepackage{color}
\usepackage{tikz}
\usetikzlibrary{arrows,shapes,decorations,decorations.text}
\usepackage{enumerate}
\usepackage{amsthm}
\usepackage{graphicx}
\usepackage{bbm}

\usepackage[shortlabels]{enumitem}

\theoremstyle{definition}
\newtheorem{theorem}{Theorem}
\numberwithin{theorem}{section}
 
 \theoremstyle{definition}
\newtheorem{lemma}[theorem]{Lemma}

\theoremstyle{definition}
\newtheorem{sublemma}{Sublemma}[theorem]

\theoremstyle{definition}
\newtheorem{definition}[theorem]{Definition}

\theoremstyle{definition}
\newtheorem{claim}[theorem]{Claim}

\theoremstyle{definition}
\newtheorem{proposition}[theorem]{Proposition}

\theoremstyle{definition}
\newtheorem{corollary}[theorem]{Corollary}

\theoremstyle{definition}

\theoremstyle{definition}

\theoremstyle{definition}
\newtheorem{remark}[theorem]{Remark}

\renewenvironment{proof}{\vskip.01in\noindent\textsc{Proof.}}{\hfill $\Box$ \vskip.1in}

\newcommand{\sym}{\rm Sym}

\newcommand{\Orb}{{\rm Orb}}
\newcommand{\dom}{{\rm dom}}

\newcommand{\setm}{\mathop{\setminus}}
\newcommand{\fix}{{\rm fix}}
\newcommand{\Sym}{{\rm Sym}}

\newcommand{\rng}{{\rm range}}

\newcommand{\DM}{\mathsf{EDM_{DM}}}
\newcommand{\BG}{\mathsf{EDM_{BG}}}
\newcommand{\T}{\mathsf{EDM_{T}}}
\newcommand{\bc}{{\rm bc}}

\newtheoremstyle{remark}
  {}{}{}{}{\bfseries}{.}{.5em}{{\thmname{#1 }}{\thmnumber{#2}}{\thmnote{ (#3)}}}
\newcommand{\defeq}{\stackrel{\textup{\tiny def}}{=}}

\begin{document}

\title{Three forms of the Erd\H{o}s-Dushnik-Miller Theorem}
\author[P. Howard]{Paul Howard}
\address{2770 Ember Way, Ann Arbor, MI 48104, USA}
\email{phoward@emich.edu}
\author[E. Tachtsis]{Eleftherios Tachtsis}
\address{Department of Statistics and Actuarial-Financial Mathematics, University of the Aegean, Karlovassi 83200, Samos, Greece}
\email{ltah@aegean.gr}

\begin{abstract}
We continue the study of the Erd\H{o}s-Dushnik-Miller theorem (A graph with an uncountable set of vertices has either an infinite independent set or an uncountable clique) 
in set theory without the axiom of choice.  We show that there are three inequivalent versions of this theorem and we give some results about the positions of these versions in the deductive hierarchy of weak choice principles. Furthermore, we settle some open problems from Tachtsis \cite{Tachtsis-2024a} and from Banerjee and Gopaulsingh \cite{Banerjee-Gopaulsingh-2023}.
\end{abstract}



\maketitle

\tableofcontents

\section{Introduction} \label{S:Intro}
The purpose of this paper is to continue the study of the deductive strength of the Erd\H{o}s-Dushnik-Miller Theorem in set theory without the axiom of choice.  The recent papers of Tachtsis \cite{Tachtsis-2022b} and \cite{Tachtsis-2024a} and of Banerjee and Gopaulsingh \cite{Banerjee-Gopaulsingh-2023} have made major strides in this area.
\par
There are several equivalent ways of stating the theorem.  For example, Dushnik and Miller's Theorem 5.23 in \cite{Dushnik-Miller-1941} is
\begin{theorem}[$\mathsf{EDM}$]
\label{T:EDM}
Any infinite graph $G = (V,E)$ not containing an independent set of size $\aleph_0$ contains a complete subgraph of size $|V|$.
\end{theorem}
(Definitions will be given in the next section.)
\par 
Since Dushnik and Miller were working in \textsf{ZFC} (Zermelo-Fraenkel set theory with the axiom of choice, \textsf{AC}), the various possible definitions of ``infinite'', which may be inequivalent in the absence of \textsf{AC}, were not considered.  We also note that Theorem \ref{T:EDM} for the case that $V$ is countable had been proved earlier by Ramsey in \cite{Ramsey-1929}.  We therefore follow the convention of \cite{Tachtsis-2024a} and \cite{Banerjee-Gopaulsingh-2023} and consider the theorem only in the case where $V$ is uncountable.

So when studying the strength of \textsf{EDM} in \textsf{ZF} (\textsf{ZFC} without choice), Tachtsis \cite{Tachtsis-2024a} uses the following form:
\par 
\vskip.1in
\noindent $\mathsf{EDM}( {\nleq} \aleph_0, 
{\nless} \aleph_0, {\nleq} \aleph_0)$: If $G =(V,E)$ is a graph with an uncountable set of vertices (that is, $|V| \nleq \aleph_0$) then either $V$ contains an infinite, independent set $I$ of vertices (that is, $|I| \nless \aleph_0$ and no pair of distinct elements of $I$ is in $E$) or there is a subgraph $G' = (V', E')$ of $G$ such that $V'$ is uncountable and $G'$ is a clique ($|V'| \nleq \aleph_0$ and every pair of distinct vertices from $V'$ is in $E'$).
\vskip.1in
In the notation we have adopted, the first argument of \textsf{EDM} describes the size of $V$ (either ${\nleq} \aleph_0$ or ${>} \aleph_0$), the second argument is the size of $I$ (either ${\nless} \aleph_0$ or ${\geq} \aleph_0$) and the third argument is the size of $V'$ (either ${\nleq} \aleph_0$ or ${>} \aleph_0$ or ${=}|V|$).  Using this notation, the original theorem as it appears in \cite{Dushnik-Miller-1941} (omitting the countable case) could be interpreted as
\par 
\vskip.1in
\noindent $\mathsf{EDM}({\nleq} \aleph_0, {\ge} \aleph_0, {=} |V|)$:  If $G = (V,E)$ is a graph such that $|V| \nleq \aleph_0$ then either $V$ contains an independent $I$ such that $|I| \ge \aleph_0$ or there is a subgraph $G' = (V', E')$ of $G$ such that $|V'| = |V|$ and $G'$ is a clique.
\vskip.1in
\noindent and the version appearing in \cite{Banerjee-Gopaulsingh-2023} is $\mathsf{EDM}({\nleq} \aleph_0, {\geq}\aleph_0, {\nleq} \aleph_0)$. 
\par 
Beginning at the end of Section \ref{S:Reducing} we will refer to these three versions of $\mathsf{EDM}$ as $\T$, $\DM$ and $\BG$ respectively.  (The reasons for this are given in Section \ref{S:Reducing}.)  In this paper we study the position of these and other versions of \textsf{EDM} in the deductive hierarchy of weak choice principles in \textsf{ZF}.

\section{Definitions}

\begin{definition} \label{D:mathscrP}
Assume $X$ is a set.  Following the notation of Jech \cite{Jech} we let
\begin{enumerate}
\item $\mathscr{P}(X)$ be the power set of $X$ ($= \{y : y \subseteq X \}$) and
\item for $\alpha$ an ordinal, $\mathscr{P}^\alpha(X)$ is defined by 
  \begin{enumerate}
  \item $\mathscr{P}^0(X) = X$, 
  \item $\mathscr{P}^{\alpha}(X) = \mathscr{P}^\beta(X) \cup \mathscr{P}(\mathscr{P}^\beta(X))$ if $\alpha = \beta + 1$,
  \item $\mathscr{P}^\alpha(X) = \bigcup_{\beta < \alpha} \mathscr{P}^\beta (X)$ if $\alpha$ is a limit ordinal.
  \end{enumerate}
  \item \label{DP:mathscr4} $\mathscr{P}^\infty(X) = \displaystyle{\bigcup_{\alpha \in On} \mathscr{P}^\alpha(X)}$ where $On$ is the class of ordinals.
\end{enumerate}
\end{definition}
\begin{definition}
A set $X$ is called:
\begin{enumerate}
\item denumerable or countably infinite if $|X| = \aleph_0$;
\item countable if $|X| \le \aleph_0$;
\item uncountable if $|X| \nleq \aleph_0$.
\item infinite if $|X| \nless \aleph_0$;
\item Dedekind infinite if $|X| \ge \aleph_0$.
\end{enumerate}
\end{definition}
\begin{definition} \label{D:GraphStuff}
Let $G = (V,E)$ be a graph (that is, $V$ is a set called the set of vertices of $G$ and $E$ is a set of unordered pairs $\{v_1, v_2 \}, v_1 \ne v_2$ of elements of $V$ called the set of edges of $G$.)
\begin{enumerate}
\item $u$ and $v$ in $V$ are \emph{adjacent} if $\{u, v\} \in E$.
\item A set $W \subseteq V$ is \emph{independent} or an \emph{anticlique} if for all $u$ and $v$ in $W$, $\{u,v\} \notin E$.
\item $G$ is a \emph{complete graph} or a \emph{clique} if any two different vertices in $G$ are adjacent.
\item $H = (W,F)$ is a \emph{subgraph} of $G$ if $W \subseteq V$ and $F \subseteq E$.
\end{enumerate}

\end{definition}

\section{Weak Forms of $\mathsf{AC}$}
In the list below, for those weak forms of $\mathsf{AC}$ that appear in Howard--Rubin \cite{Howard-Rubin-1998}, we provide their form number in \cite{Howard-Rubin-1998} to facilitate the interested reader in obtaining further information on their deductive strength using \cite{Howard-Rubin-1998}.   
\begin{enumerate}
\item (Form 9) $\mathsf{DF {=} F}$:  Every Dedekind finite set $X$ ($|X| \ngeq \aleph_0$) is finite ($|X| < \aleph_0$).

\item (Form 202) $\mathsf{AC^{LO}}$: Every linearly orderable set of
non-empty sets has a choice function.

\item (Form 10) $\mathsf{AC^{\aleph_0}_{fin}}$: Every countable set of non-empty finite sets has a choice function.

\item (Form [10 E]) $\mathsf{PAC^{\aleph_0}_{fin}}$:  Every countably infinite set of non-empty finite sets has an infinite subset with a choice function.\footnote{Such a function is called a \emph{partial choice function} for the original (infinite) set.}


\item Let $\kappa$ be a well-ordered cardinal with $\kappa\ge 2$ and let $\lambda$ be an aleph. $\mathsf{AC}_{\kappa}^{\lambda}$: Every $\lambda$-sized set of sets each of size $\kappa$ has a choice function. In particular, $\mathsf{AC}^{\aleph_0}_{\aleph_0}$ is Form [32 A] in \cite{Howard-Rubin-1998}.

\item (Form [32 B]) $\mathsf{PAC}^{\aleph_0}_{\aleph_0}$: Every countably infinite set of countably infinite sets has an infinite subset with a choice function.

\item (Form 60) $\mathsf{AC_{WO}}$:  Every family of non-empty well-orderable sets has a choice function.

\item (Form [18 A]) $\mathsf{PC}(\aleph_0,2,\aleph_0)$:  Every countably infinite set of $2$-element sets has an infinite subset with a choice function.

\item (Form 373($n$)) For $n\in\omega$, $n\ge 2$, $\mathsf{PC}(\aleph_0,n,\aleph_0)$:  Every countably infinite set of $n$-element sets has an infinite subset with a choice function.

\item $\mathsf{PC}({>} \aleph_0, {<} \aleph_0, {>} \aleph_0)$:  Every family $X$ of non-empty finite sets such that $|X| > \aleph_0$ has a subfamily $Y$ such that $|Y| > \aleph_0$ and $Y$ has a choice function.

\item $\mathsf{PC}({\nleq} \aleph_0, {<} \aleph_0, {>} \aleph_0)$:  Every family $X$ of non-empty finite sets such that $|X| \nleq \aleph_0$ has a subfamily $Y$ such that $|Y| > \aleph_0$ and $Y$ has a choice function.

\item $\mathsf{PC}({>} \aleph_0, {<} \aleph_0, {\nleq} \aleph_0)$:  Every family $X$ of non-empty finite sets such that $|X| > \aleph_0$ has a subfamily $Y$ such that $|Y| \nleq \aleph_0$ and $Y$ has a choice function.

\item $\mathsf{PC}(\not\leq\aleph_{0},<\aleph_{0},\not\leq\aleph_{0})$: Every family $X$ of non-empty finite sets such that $|X| \nleq \aleph_0$ has a subfamily $Y$ such that $|Y| \nleq \aleph_0$ and $Y$ has a choice function.

\item Kurepa's Theorem:  If $(P, \le)$ is a partially ordered set in which all anti-chains are finite and all chains are countable then $P$ is countable.

\item Let $\kappa$ be an aleph. $\mathsf{UT}(\aleph_{0},\kappa,\kappa)$: The union of a countably infinite family of sets each of size $\kappa$ has size $\kappa$.
In particular, $\mathsf{UT}(\aleph_{0},\aleph_{0},\aleph_{0})$ is the Countable Union Theorem (Form 31 in \cite{Howard-Rubin-1998}). 

\end{enumerate}

\section{Reducing the Number of Forms of $\mathsf{EDM}$ to Three} \label{S:Reducing}

Using the notation described in Section \ref{S:Intro} there are twelve possible forms of $\mathsf{EDM}$: $\mathsf{EDM}(A,B,C)$ where $A$ is ${\nleq} \aleph_0$ or ${>} \aleph_0$, $B$ is ${\nless}\aleph_0$ or  ${\geq} \aleph_0$ and $C$ is ${\nleq} \aleph_0$, ${>} \aleph_0$ or ${=}|V|$.  We show in this section that, under the relation ``is equivalent to in $\mathsf{ZF}$'' there are at most three equivalence classes of  these 12 forms.  We also prove several implications between forms of $\mathsf{EDM}$ and other weak forms of $\mathsf{AC}$.

\begin{proposition} \label{P:Easy1}
Let $*$ be either of ${\nless} \aleph_0$ or ${\ge} \aleph_0$ and let $**$ be any one of ${\nleq} \aleph_0$, ${>} \aleph_0$ or ${=} |V|$ then $\mathsf{EDM}({\nleq} \aleph_0, *, **)$ is equivalent to $\mathsf{EDM}({>} \aleph_0, *, **)$.
\end{proposition}
\begin{proof}
It is clear that $\mathsf{EDM}({\nleq} \aleph_0, *, **)$ implies $\mathsf{EDM}({>} \aleph_0, *, **)$.  For the other implication assume $\mathsf{EDM}({>} \aleph_0, *, **)$ and let $G  = (V,E)$ be a graph for which $|V| \nleq \aleph_0$.  Let $V' = V \cup W$ where $W$ is a set  for which $|W| = \aleph_0$ and $V \cap W = \emptyset$.  Let $E' = E \cup \{ \{ w_1, w_2 \} : w_1, w_2 \in W \land w_1 \ne w_2 \}$ and let $G' = (V', G')$.  By $\mathsf{EDM}({>} \aleph_0, *, **)$, either 
\begin{enumerate}
\item \label{I:I} $V'$ has a subset $I'$ such that $*$ is true of $I'$ and $I'$ is independent in $G'$ or
\item \label{I:Gprime} there is a subgraph $G'' = (V'', E'')$ of $G'$ such that $G''$ is a clique and $**$ is true of $|V''|$. (That is, if $**$ is ${\nleq} \aleph_0$ the $|V''| \nleq \aleph_0$, if $**$ is ${>} \aleph_0$ then $|V''| > \aleph_0$ and if $**$ is ${=} |V|$ then $|V''| = |V'|$.)
\end{enumerate}
If (\ref{I:I}) holds then $I'$ contains at most one element of $W$.  So in this case $*$ is true of $I \defeq I' \setm W$ (whether $*$ is ${\nless} \aleph_0$ or ${\ge} \aleph_0$) and $I \subseteq V$.  So the conclusion of $\mathsf{EDM}({\nleq} \aleph_0, *, **)$ is true for $G$.
\par
On the other hand, if (\ref{I:Gprime}) is true, we note that, since $G''$ is a clique, either $V'' \subseteq V$ or $V'' \subseteq W$.  But $V'' \subseteq W$ is not possible. For if $V'' \subseteq W$ then 
\begin{equation} \label{E:VleAleph0}
|V''| \le \aleph_0
\end{equation}
and this contradicts all three of the possible choices for $**$.  (If $**$ is true of $|V''|$ and $**$ is ${\leq} \aleph_0$ then $|V''| \nleq \aleph_0$.  If $**$ is ${>} \aleph_0$ then $|V''| > \aleph_0$.  And if $**$ is ${=} |V|$ then $|V''| = |V'|$ but we have assumed $|V| \nleq \aleph_0$ and therefore $|V'| =|V \cup W| \nleq \aleph_0$.)
\par 
We conclude that $V'' \subseteq V$.  It follows that $E'' \subseteq E$    so that $G''$ is a subgraph of $G$. We also know that  $V''$ is a clique.  Therefore to complete the proof we argue that $|V''|$ satisfies $**$ for each the three possible choices of $**$.
\par 
It is clear that if $**$ is either ${\nleq} \aleph_0$ or ${>} \aleph_0$  then $**$ is true of $|V''|$ (see item \ref{I:Gprime} above).  If $**$ is ${=} |V|$ then we have  $|V''| = |V'|$.  Since $V' \supset V$ this gives us $|V''| \ge |V|$.  Similarly, since $V'' \subseteq V$ we have $|V''| \leq |V|$.  So $|V''| = |V|$.  This completes the proof of the proposition.
\end{proof}
\par 
By Proposition \ref{P:Easy1} the first argument of $\mathsf{EDM}$ doesn't matter and we will omit this argument for the remainder of the paper and, when necessary, assume it is ${\nleq}\aleph_0$.
\par 
\begin{proposition} \label{P:impliesDF=F}
\begin{enumerate}
\item \label{I:idf=e1} $\mathsf{EDM}({\nless} \aleph_0, {>} \aleph_0)$ implies $\mathsf{DF{=}F}$ 
\item \label{I:idf=e2} $\mathsf{EDM}({\geq} \aleph_0, \nleq \aleph_0)$ implies $\mathsf{DF{=}F}$.
\end{enumerate}
\end{proposition}
\begin{proof}
For part (\ref{I:idf=e1}) assume $\mathsf{EDM}({\nless} \aleph_0, {>} \aleph_0)$ and let $X$ be any infinite set.  Then $|X| \nless \aleph_0$ and we have to show that $X$ is Dedekind infinite, that is, we have to show that $\aleph_0 \leq |X|$. Assume this is not the case. Then $|X| \ngeq \aleph_0$ and therefore $|X| \nleq \aleph_0$. Let $E = \{ \{x_1, x_2\} : x_1 \ne x_2 \mbox{ and } x_1, x_2  \in X \}$ then the graph $G = (X,E)$ satisfies the hypotheses of $\mathsf{EDM}({\nless} \aleph_0, {>} \aleph_0)$. Since the only independent sets in $G$ are singletons, there is a subgraph $G'  = (X', E')$ of $G$ such that $|X'| > \aleph_0$ and $G'$ is a clique.  Therefore $\aleph_0 < |X|$ contradicting our assumption.
\par 
Similarly for Part (\ref{I:idf=e2}) we assume that $X$ is infinite and not Dedekind infinite from which it follows, as in the proof of Part (\ref{I:idf=e1}), that $|X| \nleq \aleph_0$. So, letting $E = \emptyset$, the graph $G = (X,E)$ satisfies the hypotheses of $\mathsf{EDM}({\ge}\aleph_0, {\nleq} \aleph_0)$.  Since the only cliques in $G$ are singletons, there must be an independent subset $X'$ of $X$ such that $|X'| \ge \aleph_0$.  Therefore $X$ is Dedekind infinite, a contradiction.
\end{proof}
\par 
\begin{proposition} \label{P:=Vto>N0}
\begin{enumerate}
\item \label{I:=Vto>aleph} $\mathsf{EDM}({\geq} \aleph_0, {=} |V|)$ implies $\mathsf{EDM}({\geq} \aleph_0, {>} \aleph_0)$.
\item \label{I:=Vto>aleph2} $\mathsf{EDM}({\nless} \aleph_0, {=} |V|)$ implies $\mathsf{EDM}({\nless} \aleph_0, {>} \aleph_0)$.
\end{enumerate}
\end{proposition}
\begin{proof}
(\ref{I:=Vto>aleph})  Assume $\mathsf{EDM}({\geq} \aleph_0, {=} |V|)$ and let $G = (V,E)$ be a graph for which $|V| > \aleph_0$.  By $\mathsf{EDM}( {\geq} \aleph_0, {=} |V|)$ either
\begin{itemize}
\item there is an independent subset $I \subseteq V$ such that $|I| \geq \aleph_0$ in which case the conclusion of $\mathsf{EDM}({\geq} \aleph_0, {>} \aleph_0)$ is true or
\item there is a subgraph $G' = (V', E')$ of $G$ such that $G'$ is a clique and $|V'| = |V|$.  Since $|V| > \aleph_0$ we have $|V'| > \aleph_0$ and the conclusion of $\mathsf{EDM}({\geq} \aleph_0, {>} \aleph_0)$ is true in this case also.
\end{itemize}
\par 
The proof of (\ref{I:=Vto>aleph2}) is similar and we take the liberty of omitting it.
\end{proof}

\begin{proposition} \label{P:3Classes}
\begin{enumerate}
\item \label{I:=Vequiv} $\mathsf{EDM}({\nless} \aleph_0, {=}|V|)$ is equivalent to $\mathsf{EDM}({\geq} \aleph_0, {=} |V|)$.
\item \label{I:Otherequiv} $\mathsf{EDM}({\geq} \aleph_0, {\nleq} \aleph_0)$ is equivalent to each of:  $\mathsf{EDM}({\geq} \aleph_0, {>} \aleph_0)$ and $\mathsf{EDM}({\nless} \aleph_0, {>} \aleph_0)$.
\item \label{I:imps} $\mathsf{EDM}({\nless} \aleph_0, {=} |V|) \Rightarrow \mathsf{EDM}( {\nless} \aleph_0, {>} \aleph_0) \Rightarrow \mathsf{EDM}( {\nless} \aleph_0, {\nleq} \aleph_0)$.
\end{enumerate}
\end{proposition}
\begin{proof}
In the following diagram the implications represented by the top two down arrows follow from Proposition \ref{P:=Vto>N0}.  The other implications represented are clear.  For example, $\mathsf{EDM}({\geq} \aleph_0, {=} |V|)$ implies $\mathsf{EDM}({\nless} \aleph_0, {=}|V|)$ since for any set $X$, $|X| \geq \aleph_0$ implies $|X| \nless \aleph_0$.
\begin{center}
\begin{tikzpicture}
\node (Ege=) at (-2,2) {$\mathsf{EDM}({\geq} \aleph_0,{=} |V|)$};
\node (En<=) at (2,2) {$\mathsf{EDM}({\nless} \aleph_0, {=} |V|)$};
\node (Ege>) at (-2,0) {$\mathsf{EDM}({\geq} \aleph_0, {>} \aleph_0)$};
\node (En<>) at (2,0) {$\mathsf{EDM}({\nless} \aleph_0, {>} \aleph_0)$};
\node (Egenle) at (-2,-2) {$\mathsf{EDM}({\geq} \aleph_0, {\nleq} \aleph_0)$};
\node (En<nle) at (2,-2) {$\mathsf{EDM}({\nless} \aleph_0, {\nleq} \aleph_0)$};
\draw [thick, ->](Ege=) -- (En<=);
\draw [thick, ->](Ege>) -- (En<>);
\draw [thick, ->](Egenle) -- (En<nle);
\draw [thick, ->](Ege=) -- (Ege>);
\draw [thick, ->](Ege>) -- (Egenle);
\draw [thick, ->](En<=) -- (En<>);
\draw [thick, ->](En<>) -- (En<nle);
\end{tikzpicture}
\end{center}
Also, by Proposition \ref{P:impliesDF=F}, every form of $\mathsf{EDM}$ appearing in the diagram other than $\mathsf{EDM}({\nless} \aleph_0, {\nleq} \aleph_0)$ implies $\mathsf{DF {=} F}$.  Since $\mathsf{DF {=} F}$ implies that for every set $X$, $|X| \ge \aleph_0$ is equivalent to $|X| \nless \aleph_0$ and $|X| > \aleph_0$ is equivalent to $|X| \nleq \aleph_0$ we can conclude that parts (\ref{I:=Vequiv}) and (\ref{I:Otherequiv}) of the proposition are true.  Part (\ref{I:imps}) follows from the diagram.
\end{proof}
By Proposition \ref{P:3Classes} there are only three forms of $\mathsf{EDM}$ which might not be equivalent (and we will show below that they are not):  They are (choosing one from each equivalence class) $\mathsf{EDM}({\geq} \aleph_0, {=}|V|)$, $\mathsf{EDM}({\geq} \aleph_0, {
\nleq} \aleph_0)$ and $\mathsf{EDM}({\nless} \aleph_0, {\nleq} \aleph_0)$.  The first is the version of $\mathsf{EDM}$ proved by Dushnik and Miller and for the remainder of the paper we will use the notation $\DM$ for it.  Similarly, the second is the version studied by Banerjee and Gopaulsingh and we will use the notation $\BG$ and the third is Tachtsis' version which we will refer to as $\T$.
\par 
So
\begin{definition}
(Forms of $\mathsf{EDM}$)
\begin{enumerate}
\item $\DM$:  If $G = (V,E)$ is a graph such that $|V| \nleq \aleph_0$ then either $V$ contains an independent  set $I$ such that $|I| \ge \aleph_0$ or there is a subgraph $G' = (V', E')$ of $G$ such that $|V'| = |V|$ and $G'$ is a clique.
\item $\BG$:  If $G = (V,E)$ is a graph such that $|V| \nleq \aleph_0$ then either $V$ contains an independent  set $I$ such that $|I| \ge \aleph_0$ or there is a subgraph $G' = (V', E')$ of $G$ such that $|V'| \nleq \aleph_0$ and $G'$ is a clique.
\item $\T$:  If $G = (V,E)$ is a graph such that $|V| \nleq \aleph_0$ then either $V$ contains an independent  set $I$ such that $|I| \nless \aleph_0$ or there is a subgraph $G' = (V', E')$ of $G$ such that $|V'| \nleq \aleph_0$ and $G'$ is a clique.
\end{enumerate}
\end{definition}

\section{$\mathsf{EDM}$ and $\mathsf{PC}$: Positive Results}

\begin{proposition}
\label{p:1}
The following hold:
\begin{enumerate}
\item \label{p:1_1} $\mathsf{PC}(>\aleph_{0},<\aleph_{0},>\aleph_{0})$ is equivalent to $\mathsf{PC}(>\aleph_{0},<\aleph_{0},\not\leq\aleph_{0})\wedge\mathsf{AC_{fin}^{\aleph_{0}}}$.

\item \label{p:1_2} $\mathsf{PC}({>}\aleph_{0},{<}\aleph_{0},{\not\leq}\aleph_{0})$ is equivalent to $\mathsf{PC}({\not\leq}\aleph_{0},{<}\aleph_{0},{\not\leq}\aleph_{0})$.
\end{enumerate}
\end{proposition}

\begin{proof}
(1) ($\Rightarrow$) It suffices to show that $\mathsf{PC}(>\aleph_{0},<\aleph_{0},>\aleph_{0})$ implies $\mathsf{AC_{fin}^{\aleph_{0}}}$, since the other implication is straightforward.  And, since $\mathsf{AC_{fin}^{\aleph_{0}}}$ is equivalent to its partial version $\mathsf{PAC_{fin}^{\aleph_{0}}}$ (see \cite{Howard-Rubin-1998}), it is enough to show that  $\mathsf{PC}(>\aleph_{0},<\aleph_{0},>\aleph_{0})$ implies $\mathsf{PAC_{fin}^{\aleph_{0}}}$. 

 To this end, let $\mathcal{A}=\{A_{n}:n\in\omega\}$ be a denumerable family of pairwise disjoint non-empty finite sets. We will show that $\mathcal{A}$ has a partial choice function.  Let $\mathcal{Y} = \{ \{ a \} : a \in \bigcup \mathcal{A} \}$ and let $\mathcal{Z}  = \mathcal{Y} \cup \mathcal{A}$.
 
Clearly, $|\mathcal{Z}|\geq \aleph_{0}$. If $|\mathcal{Z}|=\aleph_{0}$, then $|\mathcal{Y}|= \aleph_{0}$, and so $\mathcal{A}$ has a choice function. Suppose that $|\mathcal{Z}|> \aleph_{0}$. By $\mathsf{PC}(>\aleph_{0},<\aleph_{0},>\aleph_{0})$, there exists $\mathcal{W}\subseteq\mathcal{Z}$ with $|\mathcal{W}|>\aleph_{0}$ and $\mathcal{W}$ has a choice function. There are two cases:
\smallskip

(i) $|\mathcal{W}\cap\mathcal{A}|=\aleph_{0}$. Then, clearly, $\mathcal{A}$ has a partial choice function.
\smallskip

(ii) $|\mathcal{W}\cap\mathcal{A}|<\aleph_{0}$. As $\mathcal{W}\subseteq\mathcal{Z}=\mathcal{Y}\cup\mathcal{A}$ and $\mathcal{Y}\cap\mathcal{A}=\emptyset$, it follows that $\mathcal{Y}$ is Dedekind-infinite so $\bigcup\mathcal{A}$ is Dedekind-infinite.  Since $\mathcal{A}$ consists of finite sets, we conclude that $\mathcal{A}$ has a partial choice function, as required.
\smallskip

($\Leftarrow$) Let $X$ be a family of pairwise disjoint non-empty finite sets with $|X|>\aleph_{0}$. Let $Y$ be a denumerable subset of $X$. By $\mathsf{AC_{fin}^{\aleph_{0}}}$, $Y$ has a choice function, say $f$. By $\mathsf{PC}(>\aleph_{0},<\aleph_{0},\not\leq\aleph_{0})$, there exists $Z\subseteq X$ with $|Z|\not\leq\aleph_{0}$ and $Z$ has a choice function, say $g$. Clearly, $|Y\cup Z|>\aleph_{0}$ and, using $f$ and $g$, we easily obtain that $Y\cup Z$ has a choice function.
\smallskip

(2) Since ``$\Leftarrow$'' is straightforward, we only argue for ``$\Rightarrow$''. Let $\mathcal{B}$ be a family of pairwise disjoint non-empty finite sets with $|\mathcal{B}|\nleq\aleph_{0}$. Let $\mathcal{V}=\{\{n\}:n\in\omega\}$ and let $\mathcal{U}=\mathcal{B}\cup\mathcal{V}$. Since $|\mathcal{B}|\nleq\aleph_{0}$ and $|\mathcal{V}|=\aleph_{0}$, $|\mathcal{U}|>\aleph_{0}$. By $\mathsf{PC}(>\aleph_{0},<\aleph_{0},\not\leq\aleph_{0})$, there exists $\mathcal{R}\subseteq\mathcal{U}$ such that $|\mathcal{R}|\nleq\aleph_{0}$ and $\mathcal{R}$ has a choice function. Clearly, $|\mathcal{R}\cap \mathcal{B}|\nleq\aleph_{0}$ and $\mathcal{R}\cap\mathcal{B}$ has a choice function, i.e. $\mathsf{PC}({\not\leq}\aleph_{0},{<}\aleph_{0},{\not\leq}\aleph_{0})$ holds, as required. 
\end{proof}

\begin{remark} \label{R:N2}
The proof of ``$\Rightarrow$'' of Proposition \ref{p:1}(\ref{p:1_1}) implicitly shows that it is relatively consistent with $\mathsf{ZF}$ that there exists an ${>}\aleph_{0}$-sized family $\mathcal{Z}$ of non-empty finite sets which has a $\not\leq\aleph_{0}$-sized subfamily $\mathcal{Y}$ with a choice function, but $\mathcal{Y}$ (and thus $\mathcal{Z}$) has no ${>}\aleph_{0}$-sized subfamily with a choice function. Indeed, this can be shown (as in Proposition \ref{p:1}) to be true in the Second Fraenkel Model and then the result can be transferred to $\mathsf{ZF}$ via the First Embedding Theorem of Jech and Sochor (see \cite[Theorem 6.1]{Jech}). See Subsection \ref{SS:N2} (``The Model $\mathcal{N}2$'') for details. 
\end{remark}

We also prove the following proposition for use in our subsection on the model $\mathcal{N}2$ (Subsection \ref{SS:N2}).

\begin{proposition} \label{P:PCnleqToPC2=}
$\mathsf{PC}({>} \aleph_0, {<} \aleph_0, {\nleq} \aleph_0)$ implies $\mathsf{PC(\aleph_0,2,\aleph_0)} \wedge \mathsf{PC(\aleph_0,4,\aleph_0)}$.
\end{proposition}

\begin{proof}
Assume  $\mathsf{PC}({>} \aleph_0, {<} \aleph_0, {\nleq} \aleph_0)$ holds. Then, by Proposition \ref{p:1}(\ref{p:1_2}), $\mathsf{PC}({\nleq} \aleph_0, {<} \aleph_0, {\nleq} \aleph_0)$ also holds.
\par
Let $\mathcal{A} = \{ A_i : i \in \omega \}$ be a denumerable family of pairwise disjoint 2 element sets.  Say $A_i = \{ a_i, b_i \}$ for each $i \in \omega$. We will show that there is an infinite subset $I$ of $\omega$ and a function $F: I \to \bigcup \mathcal{A}$ such that for all $i \in I$, $F(i) \in A_i$.
For each $i \in \omega$ we let 
$$X_i = \big{\{} \{ (a_i, a_{i+1}), (b_i, b_{i+1}) \}, \{ (a_i, b_{i+1}), (b_i, a_{i+1}) \} \big{\}} = \{ f : A_i \to A_{i+1} : f \mbox{ is one to one} \}.$$ 
Let $\mathcal{Y} = \bigcup_{i \in \omega} X_i$.  We consider two possibilities.
\vskip.1in 
\noindent \textbf{Case 1.} $|\mathcal{Y}| \leq \aleph_0$. 
\par 
Since $\mathcal{A}$ is denumerable, so is $\mathcal{Y}$. Therefore there is a well ordering $\prec$ of $\mathcal{Y}$ of type $\omega$ and we use this well ordering to obtain a choice function $g$ for $\{ X_i : i \in \omega \}$.  ($g(X_i) = $ the $\prec$-least element of $X_i$.) Fix an element of $A_0$, say $a_0$ and define a function $F : \omega \to \bigcup \mathcal{A}$ by recursion as follows:
\begin{itemize}
\item $F(0) = a_0$
\item $F(k+1) = g(X_k)(F(k))$
\end{itemize}
We can argue by induction that $F(k) \in A_k$ for all $k \in \omega$:  The fact that $F(0) \in A_0$ is clear.  Assume that $F(k) \in A_k$.  By the definition of $g$, $g(X_k) = f$ where $f \in X_k$.  That is, $f: A_k \to A_{k+1}$ and $f$ is one to one.  Therefore $f(F(k)) \in A_{k+1}$.  By the definition of $F$, 
\[
F(k+1) = g(X_k)(F(A_k)) = f(F(k)) \in A_{k+1}. 
\]
This completes the proof in Case 1.

\vskip.1in
\noindent\textbf{Case 2.} 
$|\mathcal{Y}| \nleq \aleph_0$.
\par
By $\mathsf{PC}({\nleq} \aleph_0, {<} \aleph_0, {\nleq} \aleph_0)$, we obtain a $\mathcal{W}\subseteq\mathcal{Y}$ such that $|\mathcal{W}| \nleq \aleph_0$ and $\mathcal{W}$ has a choice function, say $h$.    
Since $|\mathcal{W}|\nleq\aleph_{0}$, the set $J = \{ j \in \omega : \mathcal{W} \cap X_j \ne \emptyset \}$ is infinite and, since $|\mathcal{W} \cap X_j| \le 2$ for all $j \in \omega$, one of the two sets $J_1 = \{ j \in \omega : |\mathcal{W} \cap X_j| = 1 \}$ or $J_2 = \{ j \in \omega : |\mathcal{W} \cap X_j| = 2 \}$ is infinite.  
\vskip.1in
\noindent \textbf{Subcase 1.}  $J_1$ is infinite.
\par 
Let $I = J_1$ and for each $i \in I$, define 
\begin{enumerate}
\item  $u(i) = \bigcup ( \mathcal{W} \cap X_i )$, that is, $u(i)$ is the unique element of $\mathcal{W} \cap X_i$.
\item $F(i) = (h(u(i))_1$  where for any ordered pair $\rho$, $(\rho)_1$ is the first component of $\rho$.  
\end{enumerate}
Note that $u(i) \in X_i$ so $u(i) = f$ where $f:A_i \to A_{i+1}$ and $f$ is one to one.  Therefore $h(u(i)) = (s,t)$ for some $s \in A_i$ and some $t \in A_{i+1}$.  From this last equation it follows that $F(i) = s \in A_i$.  This completes the Subcase 1 proof.
\vskip.1in
\noindent\textbf{Subcase 2.} $J_2$ is infinite.
\par 
Assume that $j \in J_2$ and $\mathcal{W} \cap X_j = \{ f_1, f_2 \}$ where $f_1$ and $f_2$ are (the only) one to one functions from the two element set $A_j$ onto the two element set $A_{j+1}$.  Under these circumstances $f_1 \cap f_2 = \emptyset$ and therefore the relation $R_j =:\{ h(f_1), h(f_2) \} \notin \{ f_1, f_2 \}$.  Hence $\{ h(f_1), h(f_2) \}$ is not a one to one function from $A_j$ onto $A_{j+1}$.  So either the domain of the relation $R_j$ contains exactly one element (if $R_j$ is not a function) or the range of $R_j$ contains exactly one element (if $R_j$ is a function and $R_j$ is not one to one).
\par 
Since $J_2$ is infinite, one of $J_2' =: \{ j \in J_2 : |\dom(R_j)| = 1 \}$ or $J_2{''} =: \{ j \in J_2 : |\rng(R_j)| = 1 \}$ is infinite.
If $J_2'$ is infinite then we let $I = J_2'$ and for $i \in I$ we define $F(i)$ to be the unique element of $\dom(R_j)$.  Then $F(i) \in A_i$ so that $I$ and $F$ fulfill the required conditions given in the first paragraph of the proof.
\par 
If $J_2''$ is infinite let $I = \{ j+1 : j \in J_2'' \}$ and for $i \in I$, let $F(i)$ be the unique element of $\rng(R_{i-1})$.  Again we have an $I$ and an $F$ that satisfy the required conditions and the proof is complete.
\vskip.1in

For the second implication of the proposition, let $\mathcal{B}=\{B_{i}:i\in\omega\}$ be a denumerable family of pairwise disjoint $4$ element sets, say $B_{i}=\{b_{ij}:j=1,2,3,4\}$, $i\in\omega$. By the first part of the proof, $\mathsf{PC}(\aleph_0, 2, \aleph_0)$ holds. For each $i\in\omega$, we let 
$$P_{i}=\{\mathcal{R}:\mathcal{R}\text{ is a partition of $B_{i}$ into two unordered pairs}\}.$$
Note that $|P_{i}|=3$ for all $i\in\omega$. Let $\mathcal{P}=\bigcup_{i\in\omega}P_{i}$.
There are two cases.
\vskip.1in 
\noindent \textbf{Case a.} $|\mathcal{P}| \leq \aleph_0$.
\par 
Since $\mathcal{B}$ is denumerable, so is $\mathcal{P}$. Using the latter fact, we may pick, for every $i\in\omega$, an $S_{i}\in P_{i}$. Let $\mathcal{S}=\{S_{i}:i\in\omega\}$. Then $\mathcal{S}$ is a denumerable family of pairwise disjoint 2 element sets (each being a partition of $B_{i}$ into two unordered pairs). By $\mathsf{PC}(\aleph_0, 2, \aleph_0)$, there is an infinite $\mathcal{T}\subseteq\mathcal{S}$ with a choice set, say $C_{\mathcal{T}}$. Note that every member of $C_{\mathcal{T}}$ is a 2 element subset of some member of $\mathcal{B}$. Applying once more $\mathsf{PC}(\aleph_0, 2, \aleph_0)$ to the denumerable set $C_{\mathcal{T}}$, we obtain a partial choice function for $\mathcal{B}$, as required.
\vskip.1in 
\noindent \textbf{Case b.} $|\mathcal{P}| \nleq \aleph_0$.
\par
Applying $\mathsf{PC}({\nleq} \aleph_0, {<} \aleph_0, {\nleq} \aleph_0)$ to $\mathcal{P}$, we obtain a $\mathcal{V}\subseteq\mathcal{P}$ such that $|\mathcal{V}|\nleq\aleph_{0}$ and $\mathcal{V}$ has a choice function, say $f$. Since $|\mathcal{V}|\nleq\aleph_{0}$, the set $I = \{ i \in \omega : \mathcal{V} \cap P_i \ne \emptyset \}$ is infinite and, since $|\mathcal{V} \cap P_i| \le 3$ for all $i \in \omega$, one of the sets $I_k = \{ i \in \omega : |\mathcal{V} \cap P_i| = k \}$, $k=1,2,3$, is infinite.
\vskip.1in
\noindent \textbf{Subcase i.}  $I_1$ is infinite. 
\par
Let $\mathcal{N}_{1}=\{\bigcup(\mathcal{V}\cap P_{i}):i\in I_{1}\}$. Applying $\mathsf{PC}(\aleph_0, 2, \aleph_0)$ to $\mathcal{N}_{1}$, we obtain a denumerable pairwise disjoint family $\mathcal{N}_{2}$ such that every member of $\mathcal{N}_{2}$ is a 2 element subset of some $B\in\mathcal{B}$. Applying $\mathsf{PC}(\aleph_0, 2, \aleph_0)$ to $\mathcal{N}_{2}$, we obtain a partial choice function for $\mathcal{B}$, as required.
\vskip.1in
\noindent \textbf{Subcase ii.}  $I_2$ is infinite. 
\par
Let $\mathcal{M}=\{\bigcup(P_{i}\setminus\bigcup(\mathcal{V}\cap P_{i})):i\in I_{2}\}$. Since $|P_{i}|=3$ for all $i\in\omega$, it follows that, for every $i\in I_{2}$,  $\bigcup(P_{i}\setminus\bigcup(\mathcal{V}\cap P_{i}))$ is a partition of $B_{i}$ into two unordered pairs. Applying $\mathsf{PC}(\aleph_0, 2, \aleph_0)$ to $\mathcal{M}$ and then to the obtained denumerable family of pairs, gives us a partial choice function for $\mathcal{B}$ as required.
\vskip.1in
\noindent \textbf{Subcase iii.}  $I_3$ is infinite. 
\par
We describe how to construct a choice function $g$ for the denumerable subfamily $\mathcal{B}^{*}=\{B_{i}:i\in I_{3}\}$ of $\mathcal{B}$. Fix temporarily an $i\in I_{3}$. Then
$$P_{i}=\Big{\{}\big{\{}\{b_{i1}, b_{i2}\},\{b_{i3}, b_{i4}\}\big{\}}, \big{\{}\{b_{i1}, b_{i3}\},\{b_{i2}, b_{i4}\}\big{\}}, \big{\{}\{b_{i1}, b_{i4}\},\{b_{i2}, b_{i3}\}\big{\}}\Big{\}}.$$
By hypothesis, $\mathcal{V}\cap P_{i}=P_{i}$ for all $i\in I_{3}$. Moreover, via $f$ we have chosen a 2 element set from each member of $P_{i}$. Now, there are several cases how this is done. We provide 4 representative cases along with how to choose an element from $B_{i}$ in each of these cases; all  other cases (due to symmetry and the fact that every element of $P_{i}$ is a partition of $B_{i}$ into two unordered pairs) can be treated similarly and are thus left to the reader. So suppose that via $f$ we have chosen respectively from each of the above three members of $P_{i}$:
\medskip

\textbf{(a)} $\{b_{i1}, b_{i2}\}$,  $\{b_{i1}, b_{i3}\}$, $\{b_{i1}, b_{i4}\}$. Then, $\{b_{i1}, b_{i2}\} \cup \{b_{i1}, b_{i3}\} \cup \{b_{i1}, b_{i4}\}=B_{i}$ and $b_{i1}$ is the only element of $B_{i}$ which appears in each of the latter three 2 element sets. Put $g(B_{i})=b_{i1}$;
\medskip

\textbf{(b)} $\{b_{i1}, b_{i2}\}$,  $\{b_{i2}, b_{i4}\}$, $\{b_{i1}, b_{i4}\}$. Then, there is exactly one element of $B_{i}$ which is missing from all the latter three sets, namely $b_{i3}$. Put $g(B_{i})=b_{i3}$;
\medskip

\textbf{(c)} $\{b_{i1}, b_{i2}\}$,  $\{b_{i1}, b_{i3}\}$, $\{b_{i2}, b_{i3}\}$. Then, there is exactly one element of $B_{i}$ which is missing from all the latter three sets, namely $b_{i4}$. Put $g(B_{i})=b_{i4}$;
\medskip

\textbf{(d)} $\{b_{i1}, b_{i2}\}$,  $\{b_{i2}, b_{i4}\}$, $\{b_{i2}, b_{i3}\}$. Then, $\{b_{i1}, b_{i2}\} \cup \{b_{i2}, b_{i4}\} \cup \{b_{i2}, b_{i3}\}=B_{i}$ and $b_{i2}$ is the only element of $B_{i}$ which appears in each of the latter three 2 element sets. Put $g(B_{i})=b_{i2}$.
\medskip

In view of the above considerations, we may unambigiously define a function $g$ on $\mathcal{B}^{*}$ by stipulating, for all $i\in I_{3}$,
\[
g(B_{i}) = \begin{cases}  \bigcup(\bigcap f[P_{i}]) & \mbox{ if }\bigcup f[P_{i}]=B_{i}\\
                           \bigcup(B_{i}\setminus\bigcup f[P_{i}])) & \mbox{ if }\bigcup f[P_{i}] \neq B_{i} 
           \end{cases}.                           
\]
Clearly, $g$ is a choice function for $\mathcal{B}^{*}$ and, since $\mathcal{B}^{*}$ is an infinite subset of $\mathcal{B}$, $\mathcal{B}$ has a partial choice function, as required.
\par
The above arguments complete the proof of the proposition.
\end{proof}

\begin{remark}
It is worth noting that, on the basis of the proof of the second part of Proposition \ref{P:PCnleqToPC2=}, $\mathsf{PC}(\aleph_0, 2, \aleph_0) \wedge \mathsf{PC}(\aleph_0, 3, \aleph_0)$ implies $\mathsf{PC}(\aleph_0, 4, \aleph_0)$. Indeed, let $\mathcal{B}$ be a denumerable family of pairwise disjoint 4 element sets, and (as in the proof of Proposition \ref{P:PCnleqToPC2=}) let $\mathcal{P}=\{P_{B}:B\in\mathcal{B}\}$ where, for $B\in\mathcal{B}$, $P_{B}=\{\mathcal{R}:\mathcal{R}$ is a partition of $B$ into two unordered pairs$\}$. Applying first $\mathsf{PC}(\aleph_0, 3, \aleph_0)$ to $\mathcal{P}$ and then applying $\mathsf{PC}(\aleph_0, 2, \aleph_0)$ two times, readily gives us that $\mathcal{B}$ has a partial choice function, i.e. $\mathsf{PC}(\aleph_0, 4, \aleph_0)$ holds. We do not know whether or not $\mathsf{PC}({>} \aleph_0, {<} \aleph_0, {\nleq} \aleph_0)$ implies $\mathsf{PC}(\aleph_0, 3, \aleph_0)$ or $\mathsf{PC}(\aleph_0, n, \aleph_0)$ for any $n>4$. See Question \ref{Q:PC_n-element} in Section \ref{S:Questions}.
\end{remark}

\begin{proposition} \label{P:ConsequencesOfEDM}
(Consequences of the various forms of $\mathsf{EDM}$)
\begin{enumerate}
\item \label{PP:DMtoAC} $\DM$ is equivalent to $\textsf{AC}$.

\item \label{PP:BGtoPC} $\BG$ implies $\T$ which in turn implies $\mathsf{PC}({>} \aleph_0, {<} \aleph_0, {>} \aleph_0)$.

\item \label{PP:BGtoPC_2} $\BG$ implies $\mathsf{PC}({\nleq} \aleph_0, {<} \aleph_0, {>} \aleph_0)$ which in turn implies $\mathsf{PC}({>} \aleph_0, {<} \aleph_0, {>} \aleph_0)$.


\end{enumerate}
\end{proposition}

\begin{proof}
For the proof of (\ref{PP:DMtoAC}) it suffices to show that $\DM$ implies $\mathsf{AC}$ and we do this by showing that $\DM$ implies that any two sets have comparable cardinalities. (See \cite[Statement T1, page 21]{Rubin-Rubin-1985}). Let $X$ and $Y$ be disjoint sets. If both $|X| \le \aleph_0$ and $|Y| \le \aleph_0$ then $|X|$ and $|Y|$ are comparable.  Otherwise $|X \cup Y | \nleq \aleph_0$.  In this case we let $V = X \cup Y$, let $E = \{ (s,t) : s \ne t \land ( \{s, t\} \subseteq X \lor \{ s, t \} \subseteq Y) \}$ and let $G = (V,E)$.  Then $G$ satisfies the hypotheses of $\DM$.  Therefore either there is a subset $I$ of $V$ such that $|I| \ge \aleph_0$ and $I$ is independent in $G$ or $G$ has a subgraph $G' = (V', E')$ which is a clique and $|V'| = |V| = |X \cup Y|$.  Since an independent subset in $G$ can contain at most one element of $X$ and at most one element of $Y$, the second of the two alternatives must hold.  But, if $G'$ is a clique then either $V' \subseteq X$ or $V'\subseteq Y$.  In the first case $ |X \cup Y| = |V'| \le |X|$ from which we can conclude that $|Y| \le |X|$.  Similarly, in the second case we have $|X| \le |Y|$.  This completes the proof of (\ref{PP:DMtoAC}).
\smallskip
  
For (\ref{PP:BGtoPC}), note that the first implication is straightforward.  
\par 

For the second implication of (\ref{PP:BGtoPC}), we first prove the following two lemmas.

\begin{lemma}
\label{lem:EDMT_CCC}
$\mathsf{EDM_{T}}$ implies $\mathsf{AC}_{\aleph_{0}}^{\aleph_{0}}$.
\end{lemma}

\begin{proof} 
It is known that $\mathsf{AC}_{\aleph_{0}}^{\aleph_{0}}$ is equivalent to its partial version $\mathsf{PAC}^{\aleph_0}_{\aleph_0}$ (see \cite{Howard-Rubin-1998}). So, assuming $\mathsf{EDM_{T}}$, we let $\mathcal{A}=\{A_{n}:n\in\omega\}$ be a denumerable family of pairwise disjoint denumerable sets and we show that $\mathcal{A}$ has a partial choice function. Assume the contrary. Consider the graph $G=(V,E)$, where $V=\bigcup\mathcal{A}$ and $E=\bigcup\{[A_{n}]^{2}:n\in\omega\}$. As $\mathcal{A}$ has no partial choice function, $|V|\not\leq\aleph_{0}$. So, $\mathsf{EDM_{T}}$ can be applied to $G$ giving us an infinite anticlique or an uncountable clique. By definition of $E$, the first possibility (i.e. infinite anticlique) yields a partial choice function for $\mathcal{A}$ which is impossible (by assumption) and the second possibility (i.e. uncountable clique) yields for some $n\in\omega$, $A_{n}$ is uncountable which is absurd. Therefore, $\mathcal{A}$ has a partial choice function, and so $\mathsf{AC}_{\aleph_{0}}^{\aleph_{0}}$ holds.
\end{proof}
\par

The proof of Lemma \ref{TtoPCnleq} below depends heavily on ideas due to Tachtsis, see \cite[The proof of Theorem 3 that $\T$ is false in $\mathcal{N}$]{Tachtsis-2024a}.  

\begin{lemma}
\label{TtoPCnleq}
 $\T$ implies $\mathsf{PC}({>} \aleph_0, {<} \aleph_0, {\nleq} \aleph_0)$.
\end{lemma}

\begin{proof}
Assume that $X$ is a set of non-empty, pairwise disjoint, finite sets such that $|X| > \aleph_0$.  Let $G = (V,E)$ where $V = \bigcup X$ and 
\[
E = \{ \{ x, y \} \subseteq V : x \mbox{ and } y \mbox{ are in different elements of } X \}.
\]
Note that, as $X$ is disjoint and $|X|>\aleph_{0}$, we have $|V|\nleq\aleph_{0}$. Applying $\T$, since any independent set in $G$ must be finite, there is a subgraph $G' = (V',E')$ of $G$ such that $G'$ is a clique and $|V'| \nleq \aleph_0$.  Since $G'$ is a clique, $V'$ is a choice set for $Y = \{ z \in X : z \cap V' \ne \emptyset \}$ and $|Y| = |V'|$. Therefore $|Y| \nleq \aleph_0$.  
\end{proof}

\par 
By Lemma \ref{lem:EDMT_CCC}, the fact that $\mathsf{AC}_{\aleph_{0}}^{\aleph_{0}}$ implies $\mathsf{AC_{fin}^{\aleph_{0}}}$, Proposition \ref{p:1}(\ref{p:1_1}) and Lemma \ref{TtoPCnleq}, we conclude that $\mathsf{EDM_{T}}$ implies $\mathsf{PC}(>\aleph_{0},<\aleph_{0},>\aleph_{0})$, as required.

\medskip
For (3), note that the first implication follows from the fact that $\BG$ implies $\mathsf{DF=F}$ (by Proposition \ref{P:impliesDF=F}) and the fact that, under $\mathsf{DF=F}$, $|X|\nleq\aleph_{0}$ if and only if $|X|>\aleph_{0}$.

The second implication of (3) is straightforward.  
\end{proof} 

\medskip

\section{Fraenkel-Mostowski Models and $\mathsf{EDM}$}
 
  \subsection{Terminology and Properties of Fraenkel-Mostowski models}
  
{Assume that $\mathcal{M}$ is a model of $\mathsf{ZFA} + \mathsf{AC}$ whose set of atoms is $A$ and assume that $\mathcal{G}$ is a group of permutations of $A$.  For any $\phi \in \mathcal{G}$, $\phi$ can be extended to an autormorphism $\phi^*$ of $\mathcal{M}$ by $\in$-induction.  (That is, $\phi^*$ is defined by $\phi^*(a) = \phi(a)$ for $a \in A$ and $\phi^*(x) = \{ \phi^*(y) : y \in x \}$ for $x \notin A$.)  We follow the usual convention of denoting $\phi^*$ by $\phi$ when no confusion is possible.} 
 For any $x \in \mathcal{M}$  and any subgroup $H$ of $\mathcal{G}$ we let $\fix_H(x) = \{ \phi \in H : \forall t \in x, \phi(t) = t \}$, $\Orb_H(x) = \{ \phi(x) : \phi \in H \}$ and $\Sym_H(x) = \{ \phi \in H : \phi(x) = x \}$. But note that for any set $Y$, $\Sym(Y)$ (without a subscript) denotes the set of all permutations of $Y$.  For any permutation $\phi$ (of a set $Y$), $\sup(\phi) = \{ x \in Y : \phi(x) \ne x \}$.
\par  
 A \emph{normal filter of subgroups of $\mathcal{G}$} is a collection $\Gamma$ of subgroups of $\mathcal{G}$ satisfying
\begin{enumerate}
\item $\mathcal{G} \in \Gamma$,
\item If $H \in \Gamma$ and $K$ is a subgroup of $\mathcal{G}$ for which $H \subseteq K$ then $K \in \Gamma$,
\item $\Gamma$ is closed under $\cap$,
\item For all $\phi \in \mathcal{G}$ and all $H \in \Gamma$, $\phi H \phi^{-1} \in \Gamma$ and
\item $\forall a \in A$, $\Sym_{\mathcal{G}}(a) \in \Gamma$.
\end{enumerate}
Assuming that $\Gamma$ is a normal filter of subgroups of $\mathcal{G}$ and $x \in \mathcal{M}$, we say that \emph{$x$ is $\Gamma$-symmetric} if $\Sym_{\mathcal{G}}(x) \in \Gamma$ and we say that \emph{$x$ is hereditarily $\Gamma$-symmetric} if $x$ and every element of the transitive closure of $x$ is $\Gamma$-symmetric.
\par
The Fraenkel-Mostowski model determined by $\mathcal{M}$, $\mathcal{G}$ and $\Gamma$ is the class of all hereditarily $\Gamma$-symmetric sets in $\mathcal{M}$.  
\par
We refer the reader to \cite[Section 4.2]{Jech} for a more complete description of Fraenkel-Mostowski models and their properties.  In particular we will be using
\begin{lemma}
\label{L:WOSetsInFMModels}
If $\mathcal{N}$ is the Fraenkel-Mostowski model determined by $\mathcal{M}$, $G$ and $\Gamma$ then for all $x \in \mathcal{N}$ 
\begin{enumerate}
\item \label{L:FMModsPt1} $\fix(x) \in \Gamma$ if and only if $x$ is {well-orderable} in $\mathcal{N}$.
\item \label{L:FMModsPt2}If $x$ is well orderable in $\mathcal{N}$ then so is  every element of $\mathscr{P}^{\infty}(x)$. (See definition \ref{D:mathscrP} part (\ref{DP:mathscr4})).
\item \label{L:FMModsPt3} If $f \in \mathcal{M}$ is a function and $\phi \in G$ such that
  \begin{enumerate}
  \item $\phi(\dom(f)) = \dom(f)$ and
  \item $\forall z \in \dom(f)$, $\phi(f(z)) = f(\phi(z))$ 
  \end{enumerate}
then $\phi(f) = f$.  
\end{enumerate}
\end{lemma}
We will also be using the following lemma due to Banerjee and Gopaulsingh  (\cite[Proposition 3.3, (4)]{Banerjee-Gopaulsingh-2023})\footnote{In the proof in \cite{Banerjee-Gopaulsingh-2023} the hypothesis of $\BG$ (which is the same as the hypothesis of $\T$) is assumed to be true and the conclusion of $\T$ is assumed to be false.  The resulting contradiction is the negation of the conclusion of $\T$.  But the authors actually prove the negation of the conclusion of $\BG$.} and Theorem \ref{T:DNT} below.
\begin{lemma}
\label{L:BG'inFMmodels}
$\BG$ restricted to graphs with a well orderable set of vertices is true in every Fraenkel-Mostowski model.
\end{lemma}

The following theorem is due to Dixon, Neumann and Thomas \cite{Dixon-Neumann-Thomas-1986}.

\begin{theorem}
{\cite[Theorem 1]{Dixon-Neumann-Thomas-1986}} 
\label{T:DNT}
Assume $X$ is a countably infinite set and $\mathcal{L}$ is a subgroup of $\Sym(X)$ and $(\Sym(X) : \mathcal{L}) < 2^{\aleph_0}$ then there is a finite subset $S_1$ of $X$ such that $\{ \eta \in \Sym(X) : \forall x \in S_1, \eta(x) = x \} \le \mathcal{L} \le \{ \eta \in \Sym(X) : \eta[S_1] = S_1 \}$.  (Where, as usual, $\le$ is the symbol for ``is a subgroup of'' and $\eta[S_1] = \{ \eta(s) : s \in S_1 \}$.)
\end{theorem}

\subsection{The Models $\mathcal{N}1$ and $\mathcal{N}3$}
\hfill\\ \vskip-.1in                                                                                                                                                                                                                                                                                                                                                                                                                                                                                                                                                                                                                                                                                                                                                                                                                
(We will use the notation of \cite{Howard-Rubin-1998} for models that appear there.)

\par 
The Fraenkel-Mostowski model $\mathcal{N}1$ is determined by a model $\mathcal{M}$ of $\mathsf{ZFA + AC}$ with a countable set of atoms, the group $\mathcal{G} = \Sym(A)$ and the filter $\Gamma = \{ H \le \mathcal{G} : \mbox{ for some finite } E \subseteq A, \fix_{\mathcal{G}}(E) \subseteq H \}$.  This is the \emph{basic Fraenkel model}.  A description of this model can be found in \cite[Section 4.3]{Jech}.
\par 
The model $\mathcal{N}3$ is the \emph{ordered Mostowski model} described by Mostowski in \cite{Mostowski-1939} and also in \cite[Section 4.5]{Jech}.  It is determined by a model $\mathcal{M}$ of $\mathsf{ZFA + AC}$ with a countable set of atoms equipped with an ordering $\leq$ such that $(A,\leq)$ is order isomorphic to the rational numbers with the usual ordering, the group $\mathcal{G} = \{ \phi : \phi \mbox{ is an order automorphism of } (A,\leq) \}$ and the filter $\Gamma = \{ H \le \mathcal{G} : \mbox{ for some finite } E \subseteq A, \fix_{\mathcal{G}}(E) \subseteq H \}$.
\par 
\begin{theorem}
In both $\mathcal{N}1$ and $\mathcal{N}3$
\begin{enumerate}
\item $\T$ is true.
\item $\mathsf{PC}({>}\aleph_{0},{<}\aleph_{0},{>}\aleph_{0})$ is true.
\item $\mathsf{PC}({\nleq}\aleph_{0},{<}\aleph_{0},{>}\aleph_{0})$ is false.
\item $\DM$ and $\BG$ are false.
\end{enumerate}
\end{theorem}
\begin{proof}
(1) The facts that $\T$ is true in both $\mathcal{N}1$ and $\mathcal{N}3$ are proved by Banerjee and Gopaulsingh in \cite[the proofs of Theorem 4.2, parts (1) and (3)]{Banerjee-Gopaulsingh-2023}.\footnote{The authors derive a contradiction by assuming the negation of $\T$ rather than the negation of $\BG$.  In this case the proof cannot be modified to obtain a contradiction from the negation of $\BG$.}
\par
(2) This follows from (1) and Proposition \ref{P:ConsequencesOfEDM}(\ref{PP:BGtoPC}).
\par
(3) We first note that in each of $\mathcal{N}1$ and $\mathcal{N}3$, the power set, and thus the set of all finite subsets, of the set of atoms is Dedekind finite (see \cite{Jech}). This readily yields $\mathsf{PC}({\nleq}\aleph_{0},{<}\aleph_{0},{>}\aleph_{0})$ is false in $\mathcal{N}1$ and $\mathcal{N}3$. 
\par
(4) By Proposition \ref{P:ConsequencesOfEDM}(\ref{PP:DMtoAC}), $\DM$ is equivalent to $\mathsf{AC}$, and thus $\DM$ is false in $\mathcal{N}1$ and $\mathcal{N}3$. On the other hand, the facts that $\BG$ is false in $\mathcal{N}1$ and $\mathcal{N}3$ follow from part (3) and Proposition \ref{P:ConsequencesOfEDM}(\ref{PP:BGtoPC_2}).
\end{proof}

\subsection{The Model $\mathcal{N}2$} \label{SS:N2} \hfil\\ \vskip-.1in
The model $\mathcal{N}2$ is known as the \emph{Second Fraenkel Model} and is one of the models described in \cite{Howard-Rubin-1998}.  Our notation will vary slightly from the notation used in \cite{Howard-Rubin-1998}.  The ground model is a model of $\mathsf{ZFA + AC}$ with a set of atoms written as a disjoint union of pairs $A = \bigcup_{k \in \omega} B_k$ where $B_k = \{ a_k, b_k \}$ for each $k \in \omega$.  The group $\mathcal{G}$ is the group of permutations of $A$ that fix $B$ pointwise and the filter $\Gamma = \{ H \le \mathcal{G} : \mbox{ for some finite } E \subseteq A, \fix_{\mathcal{G}}(E) \subseteq H \}$.  We include $\mathcal{N}2$ because it provides some information about the relationship between $\mathsf{PC({>}\aleph_0, {<}\aleph_0, {>}\aleph_0)}$ and $\mathsf{PC({>}\aleph_0, {<}\aleph_0, {\nleq}\aleph_0)}$.  (See Proposition \ref{p:1}(\ref{p:1_1}) and Remark \ref{R:N2}.)
\par 
As described in Remark \ref{R:N2}, $\mathcal{N}2$ witnesses the fact that the following statement is not provable in $\mathsf{ZFA}$ and the result is transferable to $\mathsf{ZF}$:
\begin{center}
Every ${>}\aleph_0$-sized family $\mathcal{Z}$ of finite sets which has a ${\nleq} \aleph_0$-sized subfamily with a choice function has a ${>}\aleph_0$-sized subfamily with a choice function.
\end{center}
In $\mathcal{N}2$ we can take $\mathcal{Z} = \{B_k : j \in \omega \} \cup \{ \{ a \} : a \in A \}$.  Then, as in the proof of ``$\Rightarrow$'' in Proposition \ref{p:1}(\ref{p:1_1}), $\mathcal{Z}$ satisfies the required conditions.
\par 
However, $\mathcal{N}2$ is not a model of $\mathsf{PC}({>}\aleph_0, {<} \aleph_0, {\nleq} \aleph_0) \land \mathsf{PC}({>}\aleph_0, {<} \aleph_0, {>} \aleph_0)$ since $\mathsf{PC}(\aleph_0, 2, \aleph_0)$ is false in $\mathcal{N}2$ (see \cite{Howard-Rubin-1998}) so, by Proposition \ref{P:PCnleqToPC2=}, $\mathsf{PC}({>}\aleph_0, {<} \aleph_0, {\nleq} \aleph_0)$ is also false.  We do not know whether or not $\mathsf{PC}({>}\aleph_0, {<} \aleph_0, {\nleq} \aleph_0)$ implies $\mathsf{PC}({>}\aleph_0, {<} \aleph_0, {>} \aleph_0)$.  See Question \ref{Q:PC1toPC2} in Section \ref{S:Questions}. 

\subsection{The Model $\mathcal{N}5$ }\hfill\\ \vskip-.1in
We will use this model for a single result, namely that $\mathsf{PC}({>} \aleph_0, {<} \aleph_0, {>} \aleph_0)$ does not imply $\T$ in $\mathsf{ZF}$.  We refer the reader to \cite{Howard-Rubin-1998} for a description of $\mathcal{N}5$.  

\begin{proposition}
\label{p:3a}
$\mathsf{PC}(>\aleph_{0},<\aleph_{0},>\aleph_{0})$ (and thus $\mathsf{PC}(>\aleph_{0},<\aleph_{0},\not\leq\aleph_{0})$) is strictly weaker than $\mathsf{EDM_{T}}$ (and thus $\mathsf{EDM_{BG}}$) in $\mathsf{ZF}$.
\end{proposition}

\begin{proof}
In \cite[Theorem 6]{Tachtsis-2022b}, it was shown that $\mathsf{AC_{WO}}$  does not imply Kurepa's Theorem in $\mathsf{ZFA}$.  The model $\mathcal{N}5$ was used. Since $\mathsf{EDM_{T}}$ implies Kurepa's Theorem (see \cite{Banerjee-Gopaulsingh-2023}), it follows that $\mathsf{AC_{WO}}$ does not imply $\mathsf{EDM_{T}}$ in $\mathsf{ZFA}$. Since $\neg\mathsf{EDM_{T}}$ is a boundable statement, and thus surjectively boundable (see \cite[Note 103]{Howard-Rubin-1998} for the definition of the latter terms), and $\mathsf{AC_{WO}}$ is a class 2 statement (see \cite[p. 286]{Howard-Rubin-1998}), it follows from \cite[Theorem, top of p. 286]{Howard-Rubin-1998} that $\mathsf{AC_{WO}}\wedge\neg\mathsf{EDM_{T}}$ has a $\mathsf{ZF}$-model. As $\mathsf{AC_{WO}}$ implies $\mathsf{PC}(>\aleph_{0},<\aleph_{0},>\aleph_{0})$, $\mathsf{PC}(>\aleph_{0},<\aleph_{0},>\aleph_{0})\wedge\neg\mathsf{EDM_{T}}$ has a $\mathsf{ZF}$-model, as required.  
\end{proof}

\subsection{The Model $\mathcal{N}_T$}\hfill\\ \vskip-.1in                                                                                                                                                                                                                                                                                                                                                                                                                                                                                                                                                                                                                                                                                                                                                                                                                
\par 
We begin by describing the model.  It does not appear in \cite{Howard-Rubin-1998} but is a new model constructed by Tachtsis in \cite{Tachtsis-2024a}.
\par
The ground model $\mathcal{M}$ is a model of $\mathsf{ZFA + AC}$ with set $A$ of atoms such that $A$ is a disjoint union of pairs $A = \bigcup \{ A_i : i \in \aleph_1 \}$.  The group $\mathcal{G}$ is the group of permutations $\phi$ of $A$ such that for every $i \in \aleph_1$, there is a $j \in \aleph_1$ such that $\phi(A_i) = A_j$.  The filter $\mathcal{F}$ is the filter consisting of all subgroups $H$ of $\mathcal{G}$ such that for some countable $E \subseteq A$, $\fix_{\mathcal{G}}(E) \subseteq H$.  Tachtsis shows in \cite[Theorem 3]{Tachtsis-2024a} that $\T$ is false in this model and $\mathsf{AC^{LO}}$ is true.
\par 
Since $\mathsf{AC^{LO}}$ implies $\mathsf{DF{=}F}$ (see \cite{Howard-Rubin-1998}), $\mathsf{DF{=}F}$ is also true in the model.  We use this fact in the diagram appearing in Section \ref{S:Diagram}.

\subsection{The Model $\mathcal{N}12(\aleph_1)$}\label{SS:N12} \hfill\\ \vskip-.1in                                                                                                                                                                                                                                                                                                                                                                                                                                                                                                                                                                                                                                                                                                                                                                                                                

In this subsection we answer in the negative Question 3, part 3 in \cite{Tachtsis-2024a} by showing that $\BG$ (and therefore $\T$) is true in the model $\mathcal{N}12(\aleph_1)$.
\par
We first describe the model $\mathcal{N}12(\aleph_1)$ which appears in \cite{Howard-Rubin-1998}.
We begin with a model $\mathcal{M}$ of $\mathsf{ZFA + AC}$ which has a set of atoms $A$ of cardinality $\aleph_1$.  We let $\mathcal{G}$ be the group of all permutations of $A$ and we let $\mathcal{I}$, the ideal of supports, be the ideal of all countable subsets of $A$ so that  
\[
\mathcal{F} = \{ H \leq \mathcal{G} : \exists E \subseteq A \mbox{ such that } |E| \le \aleph_0 \land \fix_{\mathcal{G}}(E) \subseteq H \}.  
\]
The permutation model $\mathcal{N}12(\aleph_1)$ is the model determined by $\mathcal{M}$, $A$ and $\mathcal{F}$.

We now show that $\BG$ is true in $\mathcal{N}12(\aleph_1)$ beginning with some preliminary lemmas.

\begin{lemma} \label{L:countableSup}
Assume that $W$ is a countable subset of $A$ and $\tau \in \mathcal{G}$ then there is an $\alpha \in \mathcal{G}$ such that $\alpha \restriction W = \tau \restriction W$ and $\sup(\alpha)$ is countable.
\end{lemma}
\begin{proof}
Let $X = \{ \tau^n(a) : n \in \mathbb{Z} \mbox{ and } a \in W \}$ then $X$ is countable since $W$ is countable, $W \subseteq X$ and $\tau \restriction X$ is a permutation of $X$.
\par 
Define $\alpha$ by
\[
\alpha(a) = \begin{cases}  \tau(a) & \mbox{ if } a \in X \\
                           a & \mbox{ if } a \in A \setm X 
           \end{cases}                           
\]
then $\alpha \in \mathcal{G}$ and, since  $\sup(\alpha) \subseteq X$,  $\sup(\alpha)$ is countable.  Lastly it is clear that $\alpha$ agrees with $\tau$ on $W$ since $W \subseteq X$.
\end{proof}
\par 
\begin{lemma} \label{L:supEtaInD}
Assume that 
\begin{enumerate}
\item \label{I:Uprops} $U \subseteq A$ and $|U| \le \aleph_0$,
\item \label{I:tProps} $t \in \mathcal{N}12(\aleph_1)$ and $U$ is a support of $t$,
\item \label{I:gammaProps} $\gamma \in \mathcal{G}$ and $\gamma(t) \ne t$ and
\item \label{I:Dprops} $D \subseteq A \setm U$ and $|D| = \aleph_0$.
\end{enumerate}
Then there is a $\beta \in \mathcal{G}$ such that
\begin{enumerate}[(a)]
\item \label{I:supBeta} $\sup(\beta) \subseteq (\sup(\gamma) \cap U) \cup D$ and
\item \label{I:betaNet} $\beta(t) \ne t$.
\end{enumerate}
\end{lemma}
\begin{proof}
By Lemma \ref{L:countableSup} (with $\tau = \gamma$) there is an $\alpha \in \mathcal{G}$ such that $\alpha \restriction U = \gamma \restriction U$ and $|\sup(\alpha)| \le \aleph_0$.  It follows from the first equality that $\sup(\alpha) \cap U = \sup(\gamma) \cap U$ 
\par 
Let $Y = \sup(\alpha) \setm U$. Since $|Y| \le \aleph_0$ and $|D| = \aleph_0$ there is a subset $Z$ of $D$ such that $|Z| = |Y|$.  Similarly, since $|A| = \aleph_1$, there is a subset $Z'$ of $A$ such that $Z' \cap (U \cup Y \cup D) = \emptyset$ and $|Z'| = |Y|$.  Let $f : Z \to Z'$ and $g : Z' \to Y$ be one to one functions onto $Z'$ and $Y$ respectively and define elements $\psi_f$ and $\psi_g$ of $\mathcal{G}$ by
\begin{equation} \label{E:defPsi_gAnd_f}
\psi_f(a) = \begin{cases} f(a) & \mbox{ if } a \in Z \\
                          f^{-1}(a) & \mbox{ if } a \in Z' \\
                          a & \mbox{ otherwise}
            \end{cases}
       \mbox{ and } \hskip.1in  
\psi_g(a) = \begin{cases} g(a) & \mbox{ if } a \in Z' \\
                          g^{-1}(a) & \mbox{ if } a \in Y \\
                          a & \mbox{ otherwise}

            \end{cases}.           
\end{equation}
Since $Z \cap Z' = \emptyset$ and $Z' \cap Y = \emptyset$, $\psi_f$ and  $\psi_g$ are well defined and are in $\mathcal{G}$.  We also note that $\psi_f^{-1} = \psi_f$ and $\psi_g^{-1} = \psi_g$.
\par 
Define $\rho \in \mathcal{G}$ by $\rho = \psi_g \, \alpha \, \psi_g$.  We will show that 
\begin{equation}
\label{E:suprho}
\sup(\rho) \subseteq (\sup(\gamma) \cap U) \cup Z'.
\end{equation}
Since $\sup(\psi_g) = Y \cup Z'$ and $\sup(\alpha) = Y \cup (\sup(\alpha) \cap U) = Y \cup (\sup(\gamma) \cap U)$ we have 
\begin{equation} \label{E:rhoSubset}
\sup(\rho) \subseteq Y \cup Z' \cup (\sup(\gamma) \cap U).  
\end{equation}
But if $a \in Y$ then $\psi_g(a) \in Z'$ so, since $Z' \cap \sup(\alpha) = \emptyset$, $\alpha(\psi_g(a)) = \psi_g(a)$.  It follows that $\rho(a) = \psi_g(\psi_g(a)) = a$.  So $Y \cap \sup(\rho) = \emptyset$.  Therefore, using (\ref{E:rhoSubset}), we obtain (\ref{E:suprho}).
\par
Define $\beta \in \mathcal{G}$ by $\beta = \psi_f \, \rho \, \psi_f$.  We argue that $\beta$ satisfies (\ref{I:supBeta}) and (\ref{I:betaNet}) of the lemma.
\par 
First note that $\sup(\psi_f) \subseteq Z \cup Z'$ so that by (\ref{E:suprho}) $\sup(\beta) \subseteq Z \cup Z' \cup (\sup(\gamma) \cap U)$.  If $a \in Z'$ then $\psi_f(a) \in Z$. Since $Z \cap \sup(\rho) = \emptyset$, $\rho(\psi_f(a)) = \psi_f(a)$ so $\beta(a) = \psi_f(\rho(\psi_f(a))) = \psi_f(\psi_f(a) = a$.    Therefore, $\sup(\beta) \subseteq Z \cup (\sup(\gamma) \cap U) \subseteq D \cup (\sup(\gamma) \cap U)$.  This completes the proof of (\ref{I:supBeta}).
\par 
We prove (\ref{I:betaNet}) by contradiction:  Assume that $\beta(t) = t$  then 
\begin{equation}
\label{E:psi_fOftist}
\psi_f \, \psi_g \, \alpha \, \psi_g \, \psi_f (t) = t.
\end{equation}
Since $\psi_g$ and $\psi_f$ both fix $U$ pointwise and $U$ is a support of $t$, equation (\ref{E:psi_fOftist}) reduces to $\psi_f \, \psi_g \, \alpha \, (t) = t$.  This is equivalent to $\alpha(t) = \psi_g \, \psi_f(t)$ or $\alpha(t) = t$.  But $\alpha \restriction U = \gamma \restriction U$ and $U$ is a support of $t$ from which we conclude that $\alpha(t) = \gamma(t)$.  Therefore $\gamma(t) = t$  which contradicts assumption (\ref{I:gammaProps}) of the Lemma.
\end{proof}

\begin{theorem}
\label{T:EDMtrue}
$\BG$ is true in $\mathcal{N}12(\aleph_1)$.
\end{theorem}

\begin{proof}
To argue that $\BG$ is true in $\mathcal{N}12(\aleph_1)$ we let $G = (V,E)$ be a graph in the model with a set $V$ of vertices such that, in the model, $|V| \nleq \aleph_0$.  We have to argue that the conclusion of $\BG$ is true for $G$ in the model.  That is, we need to show that (in $\mathcal{N}12(\aleph_1)$) $G$ contains either  an independent set $I$ such that $|I| \ge \aleph_0$ or a complete subgraph $G' = (V', E')$ with $|V'| \nleq \aleph_0$. (Note that for any set $W$, $|W| \nleq \aleph_0$ is equivalent to $|W| > \aleph_0$ in $\mathcal{N}12(\aleph_1)$ since $\mathsf{DF{=}F}$ holds.) Let $S \in \mathcal{I}$ be a support of $G$.  If the set $V$ is well-orderable in $\mathcal{N}12(\aleph_1)$ then by Lemma \ref{L:BG'inFMmodels} the conclusion of $\BG$ is true. So we assume that $V$ is not well-orderable in $\mathcal{N}12(\aleph_1)$.  It follows that there is an element $v_0 \in V$ such that $S$ is not a support of $v_0 $ and therefore there is a $\tau \in \fix_{\mathcal{G}}(S)$ such that $\tau(v_0) \ne v_0$.  
Choose a set $S_0$ such that $S \cup S_0$ is a (countable) support of $v_0$ and $S \cap S_0 = \emptyset$.  By Lemma \ref{L:countableSup} (with $W = S \cup S_0$) we have $\exists \alpha_0 \in \mathcal{G}$ such that
\begin{align}
& \alpha_0 \restriction (S \cup S_0) = \tau \restriction (S \cup S_0), \label{E:alpha0-1} \\ 
& |\sup(\alpha_0) | \le \aleph_0 \mbox{ and } \label{E:alpha0-2} \\ 
& \alpha_0(v_0) \ne v_0. \label{E:alpha0-3}
\end{align}
((\ref{E:alpha0-3}) follows from the fact that $\tau(v_0) \ne v_0$ and $\alpha_0$ agrees with $\tau$ on a support of $v_0$.)
\par 
Choose a countably infinite set $T \subseteq A$ such that 
\begin{enumerate} 
 \item \label{D:Tp1} $S_0 \cup \sup(\alpha_0) \subseteq T$,
 \item \label{D:Tp2} $S \cap T = \emptyset$ and
 \item  \label{D:Tp3} $T \setm (S_0 \cup \sup(\alpha_0))$ is countably infinite.
\end{enumerate}
 \par 
 In what follows we will be using the following groups of permutations.
 
\begin{itemize}
\item $\mathcal{H} = \fix_{\mathcal{G}}(A\setm T) = \{ \psi \in \mathcal{G} : \sup(\psi) \subseteq T \}$ 
\item $\mathcal{H}' = \Sym(T) = \{ \psi \restriction T : \psi \in \mathcal{H} \}$
\item $\mathcal{K} = \fix_{\mathcal{H}} ( \{ v_0 \}) = 
\{ \psi \in \mathcal{H} : \psi(v_0) = v_0 \} = \{ \psi \in \mathcal{G} : \psi(v_0) = v_0 \land \sup(\psi) \subseteq T \}$. 
\item $\mathcal{K}' = \{ \psi \restriction T : \psi \in \mathcal{K} \}$.
\end{itemize}
 Using the usual group theoretic notation we let $\mathcal{H} \slash \mathcal{K}$ denote the set of left cosets of $\mathcal{K}$ in $\mathcal{H}$ and $(\mathcal{H}: \mathcal{K})$ denotes the index of $\mathcal{K}$ in $\mathcal{H}$.  That is, $(\mathcal{H}: \mathcal{K}) = |\mathcal{H}\slash \mathcal{K}|$.  Finally, if $\eta \in \mathcal{H}$ then $\eta \mathcal{K}$ is the left coset of $\mathcal{K}$ in $\mathcal{H}$ determined by $\eta$.
 
We note the following easy fact which we state as a Lemma for future use:
\begin{lemma}
\label{L:psiInK}
 For all $\psi \in \mathcal{G}$, $\psi \in \mathcal{K}$ if and only if $(\sup(\psi) \subseteq T$ and $\psi \restriction T \in \mathcal{K}')$.
\end{lemma}
In addition we will need Lemma \ref{L:orbGv_0InV} below for the proof of Theorem \ref{T:EDMtrue}.
\begin{lemma} 
\label{L:orbGv_0InV}
\begin{enumerate}
\item \label{LI:orbProp1} $\Orb_{\mathcal{H}}(v_0) \subseteq V$.
\item \label{LI:orbProp2} $|\Orb_{\mathcal{H}}(v_0)| = (\mathcal{H}:\mathcal{K})$.
\item \label{LI:orbProp3} $\Orb_{\mathcal{H}}(v_0)$ is well orderable in $\mathcal{N}12(\aleph_1)$.
\end{enumerate}
\end{lemma}
\begin{proof}
Item (\ref{LI:orbProp1}) is true because $\mathcal{H} \subseteq \fix_{\mathcal{G}}(S)$ and $S$ is a support of $G$.
\par
Item (\ref{LI:orbProp2}) follows from the fact that if $\eta_1 \mathcal{K}$ and $\eta_2 \mathcal{K}$ are  cosets of $\mathcal{K}$ in $\mathcal{H}$ then $\eta_1 \mathcal{K} = \eta_2 \mathcal{K}$ if and only if  $\eta_1(v_0) = \eta_2(v_0)$. 
\par 
 To argue for (\ref{LI:orbProp3}) we first note that for each $\eta \in \mathcal{H}$, $\eta(S \cup S_0)$ is a support of $\eta(v_0)$. Since $S_0 \subseteq T$ and $\eta$ fixes $A \setm T$ pointwise, $\eta(S \cup S_0) \subseteq T \cup S$.  It follows that for every $\eta \in \mathcal{H}$, $T \cup S$ is a support of $\eta(v_0)$.   Therefore $\Orb_{\mathcal{H}}(v_0)$ is well orderable in $\mathcal{N}12(\aleph_1)$. 
\end{proof}

\par 
We also note that there is a function from $\sym(T)$ onto $\mathcal{H} \slash \mathcal{K}$ given by $\sigma \mapsto \sigma' \mathcal{K}$ where $\sigma': A \to A$ is defined by 
\[
\sigma'(a) = \begin{cases} 
                               \sigma(a) & \mbox{if } a \in T \\
                               a & \mbox{otherwise}
                             \end{cases}.
\]                             
Therefore $(\mathcal{H} : \mathcal{K}) \le |\sym(T)| = 2^{\aleph_0}$                          
\par
We first consider the case where $(\mathcal{H}: \mathcal{K}) = 2^{\aleph_0}$.  In this case $|\Orb_{\mathcal{H}}(v_0)| = 2^{\aleph_0}$ (by Lemma \ref{L:orbGv_0InV}, part (\ref{LI:orbProp2})) and $\Orb_{\mathcal{H}}(v_0)$ is well orderable in $\mathcal{N}12(\aleph_1)$ (using the same Lemma, part (\ref{LI:orbProp3})). 
So, if we let $V' = \Orb_{\mathcal{H}}(v_0)$ and $E' = \{ \{v, w \} : v, w \in V' \land \{v, w\} \in E \}$, then the graph $G' = (V',E')$ is a subgraph of $G$, by Lemma \ref{L:orbGv_0InV}, part \ref{LI:orbProp1}, which has either a countable set of independent vertices in $\mathcal{N}12(\aleph_1)$ or an uncountable clique in $\mathcal{N}12(\aleph_1)$ using  Lemma \ref{L:BG'inFMmodels}.  It follows from the definition of $E'$ that an independent set in $G'$ is also independent in $G$ and a clique in $G'$ is also a clique in $G$.  So the graph $G$ has (in $\mathcal{N}12(\aleph_1)$) either a countably infinite set of independent vertices or an uncountable clique.  This completes the proof in our first case.
\par 
In the other possible case, where $(\mathcal{H} : \mathcal{K}) < 2^{\aleph_0}$, we use Theorem \ref{T:DNT} as follows:
For each $\psi \in \mathcal{H}$, let $I(\psi) = \psi \restriction T$.  Then $I$ is an isomorphism from $\mathcal{H}$ onto $\Sym(T) = \mathcal{H}'$.  In addition $I[\mathcal{K}] = \mathcal{K}'$.  It follows that $(\mathcal{H}':\mathcal{K}') = (\mathcal{H} : \mathcal{K}) < 2^{\aleph_0}$.  
\par
We apply Theorem \ref{T:DNT} with $X = T$ and $\mathcal{L} = \mathcal{K}'$.  This gives us a finite subset $S_1$ of $T$ such that 
\begin{equation} \label{E:S_1property}
\{ \eta \in \Sym(T): \forall s \in S_1, \eta(s) = s \} \le \mathcal{K}' \le \{ \eta \in \Sym(T) : \eta[S_1] = S_1 \}. 
\end{equation}

\begin{claim}
$S \cup S_1$ is a support of $v_0$.
\end{claim}

\begin{proof} (of the claim)
Assume $\psi \in \fix_\mathcal{G}(S \cup S_1)$ and, toward a proof by contradiction, assume that $\psi(v_0) \ne v_0$. Applying Lemma \ref{L:supEtaInD} with $U = S \cup S_0 \cup S_1$, $t = v_0$, $\gamma = \psi$ and $D = T \setm (S_0 \cup S_1)$ we get a permutation $\beta \in \mathcal{G}$ such that 
\begin{align}
 \sup(\beta) &\subseteq (\sup(\psi) \cap (S \cup S_0 \cup S_1)) \cup (T \setm (S_0 \cup S_1)) \label{E:supBetaSub}  \\
\beta(v_0) &\ne v_0. \label{E:betav_0nev_0}
\end{align}
Since $\sup(\psi) \cap (S \cup S_1) = \emptyset$ equation (\ref{E:supBetaSub}) is equivalent to
\begin{equation}
\label{E:supBetaAgain}
 \sup(\beta) \subseteq (\sup(\psi) \cap  S_0 ) \cup (T \setm (S_0 \cup S_1))
\end{equation}
Since $S_0  \subseteq T$ we conclude from (\ref{E:supBetaAgain}) that $\sup(\beta) \subseteq T$.  So
\begin{equation}
\label{E:betaRstT}
\beta \restriction T \in \Sym(T).
\end{equation}
Similarly, $\sup(\beta) \cap S_1 = \emptyset$ so $\beta \in \fix_{\mathcal{G}}(S_1)$ so
\begin{equation}
\label{E:forAllaInS_1}
\forall a \in S_1, (\beta \restriction T)(a) = a.
\end{equation}
By (\ref{E:betaRstT}) and (\ref{E:forAllaInS_1}) and using (\ref{E:S_1property}) we get $\beta \restriction T \in \mathcal{K}'$ so, by Lemma \ref{L:psiInK}, $\beta \in \mathcal{K}$ which implies $\beta(v_0) = v_0$ contradicting (\ref{E:betav_0nev_0}).
\end{proof}   

Since $S$ is not a support of $v_0$ and $S \cup S_1$ is, $S_1 \ne \emptyset$.  We fix an element $s_1 \in S_1$. For each $t \in A \setm (S \cup S_1)$ we let $\beta_t$ be the transposition $(s_1,t)$ (in the group $\mathcal{G}$).  We define $V'$ by  
\begin{equation}
V' = \{ \beta_t(v_0) : t \in A \setm (S \cup S_1) \} 
\end{equation}

\begin{claim}
\label{C:phiBeta=BetaPhi}
For all $t \in A \setm (S \cup S_1)$ and for all $\phi \in \fix_{\mathcal{G}}(S \cup S_1)$,
\[
\phi(\beta_t(v_0)) = \beta_{\phi(t)}(v_0).
\]
\end{claim}
\begin{proof}
We have to verify that $\phi \circ \beta_t$ and $\beta_{\phi(t)}$ agree on $S \cup S_1$ which is a support of $v_0$.  Assuming $s \in S \cup S_1$ there are two possibilities:  If $s \ne s_1$ then $\phi(\beta_t(s)) = s = \beta_{\phi(t)}(s)$.  If $s = s_1$ then $\phi(\beta_t(s_1)) = \phi(t) = \beta_{\phi(t)}(s_1)$.
\end{proof}

\begin{lemma}(Properties of $V'$)
\label{L:PropsOfV'}
\begin{enumerate}
\item \label{I:PropsOfV'1} $V' \subseteq V$.
\item \label{I:PropsOfV'2} $V' \in \mathcal{N}12(\aleph_1)$.
\item \label{I:PropsOfV'3} In $\mathcal{N}12(\aleph_1)$, $|V'| \nleq \aleph_0$.
\end{enumerate}
\end{lemma}
\begin{proof}
For (\ref{I:PropsOfV'1}) Assume $t \in A \setm (S \cup S_1)$.  We have to show that $\beta_t(v_0) \in V$.  Since $S_1 \subseteq T$, we have $S_1 \cap S = \emptyset$.  Therefore $s_1 \notin S$.  Since $t \notin S$, $\beta_t \in \fix_{\mathcal{G}}(S)$ and so, since $S$ is a support of $V$, $\beta_t(v_0) \in V$.
\par 
To prove (\ref{I:PropsOfV'2}) we show that $S \cup S_1$ is a support of $V'$.  Assume $\phi \in \fix_{\mathcal{G}}(S \cup S_1)$ and $t \in A \setm (S \cup S_1)$.  Then, by Claim \ref{C:phiBeta=BetaPhi}, $\phi(\beta_t(v_0)) =  \beta_{\phi(t)}(v_0) \in V'$ since $\phi(t) \in A \setm (S \cup S_1)$.
\par 
For (\ref{I:PropsOfV'3}) we first argue that the function $f: A \setm (S \cup S_1) \to V'$ defined by $f(t) = \beta_t(v_0)$ has support $S \cup S_1$ and is therefore in $\mathcal{N}12(\aleph_1)$.  Assume $(t,\beta_t(v_0)) \in f$ where $t \in A \setm (S \cup S_1)$ and assume $\phi \in \fix_{\mathcal{G}}(S \cup S_1)$.  Then $\phi(t,\beta_t(v_0)) = (\phi(t), \phi(\beta_t(v_0)) = (\phi(t), \beta_{\phi(t)}(v_0)) \in f$ where the last equality follows from Claim (\ref{C:phiBeta=BetaPhi}).
\par 
Secondly we argue that $f$ is one to one.  Assume that $f(t_1) = f(t_2)$ where $t_1$ and $t_2$ are in $A \setm (S \cup S_1)$. Using the definition of $f$ and our assumption we get $\beta_{t_1}(v_0) = \beta_{t_2}(v_0)$.  Toward a proof by contradiction that $f$ is one to one assume that $t_1 \ne t_2$.  Choose two elements $r_1 \ne r_2$ in $T \setm (S \cup S_1 \cup \{ t_1, t_2 \})$. (This is possible since $T$ is countably infinite, $T \cap S = \emptyset$ and $S_1$ is finite.)  Let $\eta$ be the product of two transpositions $\eta = (t_1, r_1)(t_2, r_2)$ (in the group $\mathcal{G}$).  Then $\eta(t_1) = r_1$ and $\eta(t_2) = r_2$. So, by Claim \ref{C:phiBeta=BetaPhi}, $\eta(\beta_{t_1}(v_0)) = \beta_{r_1}(v_0)$ and $\eta(\beta_{t_2}(v_0)) = \beta_{r_2}(v_0)$.  Using our assumption these two equations give us  $\beta_{r_1}(v_0) = \beta_{r_2}(v_0)$ and so $\beta_{r_2}^{-1} (\beta_{r_1}(v_0)) = v_0$.  We also have $\sup(\beta_{r_2}^{-1} \beta_{r_1}) = \{ s_1, r_1, r_2 \} \subseteq T$.  Therefore $\beta_{r_2}^{-1} \beta_{r_1} \in \mathcal{K}$.  Hence, by Lemma \ref{L:psiInK}, $(\beta_{r_2}^{-1} \beta_{r_1}) \restriction T \in \mathcal{K}'$.  By the second ``$\leq$'' in equation (\ref{E:S_1property}), $\beta_{r_2}^{-1}\beta_{r_1} [S_1] = S_1$.  But $s_1 \in S_1$ and $\beta_{r_2}^{-1} \beta_{r_1}(s_1) = r_1 \notin S_1$, a contradiction.  
\par 
Finally we note that $S \cup S_1$ is countable in $\mathcal{N}12(\aleph_1)$ and $|A| \nleq \aleph_0$ in $\mathcal{N}12(\aleph_1)$.  It follows that $|A \setm (S \cup S_1)| \nleq \aleph_0$ so that, in $\mathcal{N}12(\aleph_1)$, $|V'| \nleq \aleph_0$.
\end{proof}

 We complete the proof of Theorem \ref{T:EDMtrue} by showing that either 
\begin{align}
& \forall x, y \in V', \mbox{ if } x \ne y \mbox{ then } \{x,y\} \in E \hskip.2in \mbox{ or} \label{E:clique} \\
& \forall x, y \in V', \mbox{ if } x \ne y \mbox{ then } \{x,y\} \notin E. \label{E:independent}
\end{align}
Choose $x_1$ and $x_2$ in $V'$ such that $x_1 \ne x_2$ and let $x_1 = \beta_{t_1}(v_0)$ and $x_2 = \beta_{t_2}(v_0)$ where $t_1 \ne t_2$ and $t_1, t_2 \in A \setm (S \cup S_1)$.  We consider two cases.
\par \noindent \textbf{Case 1.}  $\{ x_1, x_2\} \in E$.  In this case we argue that for all $y_1 \ne y_2$ in $V'$, $\{ y_1, y_2\} \in E$. We first consider the subcase where $\{x_1, x_2 \} \cap \{ y_1, y_2 \} = \emptyset$.  Assume $y_1 = \beta_{r_1}(v_0)$ and $y_2 = \beta_{r_2}(v_0)$ where $r_1 \ne r_2$ and both are in $A \setm (S \cup S_1)$.  Let $\eta = (t_1, r_1)(t_2, r_2)$ (the product of transpositions).  $\eta \in \fix_{\mathcal{G}}(S \cup S_1)$ so $\eta(E) = E$ and, by Claim \ref{C:phiBeta=BetaPhi}, 
\[
\eta( \{ \beta_{t_1}(v_0), \beta_{t_2}(v_0) \}) = \{ \eta( \beta_{t_1}(v_0)), \eta(\beta_{t_2}(v_0)) \} = \{ \beta_{r_1}(v_0), \beta_{r_2}(v_0) \}
\]
which is in $E$ since $\{ \beta_{t_1}(v_0), \beta_{t_2}(v_0) \}$ is in $E$. 
\par 
In the subcase ``$\{ x_1, x_2\} \cap \{ y_1, y_2 \} \ne \emptyset$'' we can choose $z_1 \ne z_2$ in $V'$ so that $\{z_1, z_2 \} \cap \{ x_1, x_2, y_1, y_2 \} = \emptyset$.  Then the above argument with $y_1$ replaced by $z_1$ and $y_2$ replaced by $z_2$ gives us $\{ z_1, z_2 \} \in E$.  Applying the above argument again with $x_1$ replaced by $z_1$ and $x_2$ replaced by $z_2$ results in $\{ y_1, y_2 \} \in E$.
\par
 We have shown that in Case 1, (\ref{E:clique}) holds.  
\vskip.1in
\par \noindent \textbf{Case 2.}  $\{ x_1, x_2\} \notin E$.  In Case 2, equation (\ref{E:independent}) is true.  The proof is almost identical to the proof in Case 1 and we take the liberty of omitting it.
\end{proof}   

\subsection{The Model $\mathcal{N}9$}  \hfill\\ \vskip-.1in                                                                                                                                                                                                                                                                                                                                                                                                                                                                                                                                                                                                                                                                                                                                     In this subsection we answer in the negative Question 3, part 5 in \cite{Tachtsis-2024a} by showing that $\BG$ (and thus $\T$) is true in the model $\mathcal{N}9$.                                                                           
The model $\mathcal{N}9$ appears in \cite{Howard-Rubin-1998}\footnote{The description of $\mathcal{N}9$ which appears in \cite{Howard-Rubin-1998} is, at best, misleading.  We give the description from \cite{Halpern-Howard-1976}.} and appeared originally in \cite{Halpern-Howard-1976}.
\par 
We start with a model $\mathcal{M}$ of $\mathsf{ZFA}+\mathsf{AC}$ with a set $A$ of atoms which has the structure of the set 
\[
\omega^{(\omega)}=\{s:s:\omega\to\omega\wedge(\exists n\in\omega)(\forall j>n)(s_{j}=0)\}.
\]
We identify $A$ with this set to simplify the description of the group $\mathcal{G}$. 

For $s\in A$, the pseudo length of $s$ is the least natural number $k$ such that for all $\ell\ge k$, $s_{\ell}=0$.  A subset $B$ of $A$ is called \emph{bounded} there is an upper bound for the pseudo lengths of the elements of $B$.
  
\begin{definition}
We let $\mathcal{G}$ be the group of all permutations $\phi$ of $A$ such that the $\sup(\phi)$ ($=\{a \in A : \phi(a) \ne a \}$) is bounded. 
\end{definition}
 
For every $s\in A$ and every $n\in\omega$, let
\[
A_{s}^{n}=\{t\in A:(\forall j\ge n)(t_{j}=s_{j})\}.
\]
\begin{definition} \label{D:N9stuff} (Mostly from \cite{Halpern-Howard-1976})
Assume $s \in A$, $n,  m \in \omega$ and $n \le m$; then
\begin{enumerate}
\item $A_{s}^{n}$ is called \emph{the $n$-block containing $s$}.
\item For any $t \in A_{s}^{n}$, the \emph{$n$-block code of $t$} is the sequence  
\[
(t_n, t_{n+1}, t_{n+2}, \ldots) = (s_n, s_{n+1}, s_{n+2}, \ldots).
\]
The \emph{$n$-block code of $A_{s}^{n}$} is the $n$-block code of any of its elements.   We will denote the $n$-block code of an element $t\in A$ or an $n$-block $B$ by $\bc^{n}(t)$ or $\bc^{n}(B)$, respectively.
\item For any $t \in A_{s}^{n}$, the finite sequence $(t_0, t_1, t_2, \ldots, t_{n-1}) = t \upharpoonright n$ is called the \emph{$n$-location of $t$ (in $A_{s}^{n}$)}.
\end{enumerate}
\end{definition}

Note the following

\begin{enumerate}
\item  $A_{0}^n$ is the set of all elements of $A$ with pseudo length less than or equal to $n$.  (In the expression $A_{0}^n$, $0$ denotes the constant sequence all of whose terms are $0$.)
\item For $s\in A$ and $n,m\in\omega$ with $n\le m$, $A_{s}^{n}\subseteq A_{s}^{m}$.
\item \label{BsbsetB'} If $n \le m$ and $B$ is an $m$-block then the set of $n$-blocks contained in $B$ is a partition of $B$. (This follows from the previous item and the fact that any two different $n$-blocks are disjoint.)
\item Any $t \in A$ is the concatenation $(t \upharpoonright n)^{\frown} \bc^n(t)$ of the $n$-location of $t$ and the $n$-block code of $t$.
\end{enumerate}

\begin{definition}
For each $n \in \omega$, $\mathcal{G}_n$ is the subgroup of $\mathcal{G}$ consisting of all permutations $\phi \in \mathcal{G}$ such that 
\begin{enumerate}
\item \label{D:fixAn0} $\phi$ fixes $A^n_0$ pointwise,
\item \label{D:fixnBlocks} $\phi$ fixes the set of $n$-blocks, that is,  $A_{s}^{n} = A_{t}^{n}$ if and only if $A_{\phi(s)}^{n} = A_{\phi(t)}^{n}$,
\item \label{D:fixnlocations} for each $s \in A$, the $n$-location of $\phi(s)$ is the same as the $n$-location of $s$.
\end{enumerate}
\end{definition}

(Note that if $n\le m$, then $\mathcal{G}_{m}\subseteq \mathcal{G}_{n}$.) 
\par 
We let $\Gamma$ be the filter of subgroups of $\mathcal{G}$ generated by the groups $\mathcal{G}_n$, $n \in \omega$.  That is, $H \in \Gamma$ if and only if $H$ is a subgroup of $\mathcal{G}$ and there exists $n \in \omega$ such that $\mathcal{G}_n \subseteq H$.  It is shown in \cite{Halpern-Howard-1976} that $\Gamma$ is a normal filter.
 $\mathcal{N}9$ is the Fraenkel--Mostowski model of $\mathsf{ZFA}$ which is determined by $\mathcal{M}$, $\mathcal{G}$, and $\Gamma$.
\par 
Following are several definitions which are not needed for the description of $\mathcal{N}9$ but will be used in the proof that $\BG$ is true in this model (Theorem \ref{T:BGinN9}).
\begin{definition}
\label{D:Bnms-Gnm}
Assume $m$ and $n$ are in $\omega$ with $n < m$ and $s \in A$.
\begin{enumerate}
\item \label{D:B^n} $\mathcal{B}^n$ is the set of $n$-blocks.
\item \label{D:mBlockLocation} The \emph{$m$-block location of $A_s^n$ (in the $m$-block $A_s^m$)} is the sequence $(s_n, s_{n+1}, \ldots, s_{m-1})$.
\item \label{D:Bnms} $\mathcal{B}(n,m)$ is the set of $n$-blocks which are contained in the $m$-block $A_0^m$.
\item \label{D:Gnm} $\mathcal{G}(n, m)$ is the set $\{ \phi \in \mathcal{G}_n : \sup(\phi) \subseteq A_0^m \}$.  (Note that $\mathcal{G}(n,m) \notin \Gamma$.)
\end{enumerate}
\end{definition}
\begin{lemma}
\label{L:mBlockCodes}
Assume $m$ and $n$ are non-negative integers, with $n < m$ and assume $\phi \in \mathcal{G}_{m}$ and $s \in A$ then the $m$-block location of $A_{\phi(s)}^n$ in $A_{\phi(s)}^m$ is the same as the $m$-block location of $A_s^n$ in $A_s^m$.
\end{lemma}
We begin our study of $\mathcal{N}9$ with two results from \cite{Halpern-Howard-1976}.
\begin{theorem} \label{T:propsOfN9}
[\cite{Halpern-Howard-1976}] In $\mathcal{N}9$
\begin{enumerate}
\item \label{TI:2m=m} the statement $2m{=} m$ is true
\item \label{TI:Z+W=W} the law of infinite cardinal addition holds.  That is, for any two infinite sets $Z$ and $W$, if $|Z| \le |W|$ then $|Z| + |W| = |W|$.
\end{enumerate}
\end{theorem}
\begin{theorem}
\label{T:BGinN9}
In $\mathcal{N}9$, $\BG$ is true.
\end{theorem}
\begin{proof}
Let $G = (V,E)$ be a graph in $\mathcal{N}9$ for which $|V| \nleq \aleph_0$ (in $\mathcal{N}9$) and let $n_0$ be a positive integer for which $\mathcal{G}_{n_0} \subseteq \Sym_{\mathcal{G}}(G)$.  If for every $v \in V$, $\mathcal{G}_{n_0} \subseteq \Sym_{\mathcal{G}}(v)$ then $V$ is well orderable and we are done by Lemma \ref{L:BG'inFMmodels}.  Otherwise there is a $v_0 \in V$ and an $\alpha \in \mathcal{G}_{n_0}$ such that $\alpha(v_0) \ne v_0$.  Since $\alpha \in \mathcal{G}$ there is an $i_0 > n_0$  such that $\sup(\alpha) \subseteq A_0^{i_0}$.  Since $v_0 \in \mathcal{N}9$ there is a $j_0 > n_0$ such that $\mathcal{G}_{j_0} \subseteq \Sym_{\mathcal{G}}(v_0)$. 
Choose a positive integer $m_0$ greater than both $i_0$ and $j_0$. Let  $T = \mathcal{B}(n_0, m_0) \setminus \{ A_0^{n_0} \}$.  
\par 
Assume that $\phi \in \mathcal{G}(n_0,m_0)$ ($= \{ \phi \in \mathcal{G}_{n_0} : \sup(\phi) \subseteq A_0^{m_0} \}$ - See Definition \ref{D:Bnms-Gnm}, part \ref{D:Gnm}).  Recall that we also represent the extension of $\phi$ to the ground model $\mathcal{M}$ by $\phi$.  But for the remainder of this proof we will use $\phi^*$ for $\phi$ (the extension of $\phi$ to $\mathcal{M}$) restricted to $T$.  With this notational convention $\Sym(T) = \{ \phi^* : \phi \in \mathcal{G}(n_0, m_0) \}$.  Let 
\[
\mathcal{H} = \{ \phi^* : \phi \in \mathcal{G}(n_0, m_0) \}
\]
and let $\mathcal{K}$ be the subgroup of $\mathcal{H}$ defined by 
\[
\mathcal{K} = \{ \phi^* : \phi \in \mathcal{G}(n_0, m_0) \land \phi(v_0) = v_0 \}.
\]
\par 
For all $\phi \in \mathcal{G}_{m_0}$, $\phi \in \fix_\mathcal{G}(A_0^{m_0})$.  Therefore $A_0^{m_0}$ is well orderable in $\mathcal{N}9$.  The set of left cosets $\mathcal{H} \slash \mathcal{K}$ is in $\mathscr{P}^\infty (A_0^{m_0})$. So, by Lemma \ref{L:WOSetsInFMModels} part (\ref{L:FMModsPt2}), the set of left cosets $\mathcal{H} \slash \mathcal{K}$ is also well orderable in $\mathcal{N}9$.
\par 
We now consider two cases depending on the cardinality of $\mathcal{H} \slash \mathcal{K}$ in $\mathcal{N}9$.   
\vskip.05in\noindent\textbf{Case 1.}  $|\mathcal{H} \slash \mathcal{K}| > \aleph_0$.  In this case we let $V' = \{ \phi(v_0) : \phi \in \mathcal{G}(n_0,m_0) \}$.  We first show that $V' \in \mathcal{N}9$ and is well-orderable in $\mathcal{N}9$ by showing that for every element $x \in V'$, $\mathcal{G}_{m_0} \subseteq \Sym_{\mathcal{G}}(x)$:  Assume $x \in V'$, then $x = \phi(v_0)$ for some $\phi \in \mathcal{G}(n_0, m_0)$.  Assume $\eta \in \mathcal{G}_{m_0}$ then we have to show that $\eta(x) = x$.  We note that 
\begin{enumerate}
\item \label{I:1} For all $t \in A_0^{m_0}$, $\phi^{-1} \eta \phi (t) = t$. (Since $\eta$ is the identity on $A_0^{m_0}$.) 
\item \label{I:2} If $s \in A$ and $s$ is not the zero sequence then $\phi^{-1} \eta \phi (A_s^{m_0})$ is an $m_0$ block and for all $t \in A_s^{m_0|}$ the $m_0$ location of $\phi^{-1} \eta \phi (t)$ is the same as the $m_0$ location of $t$. (Because $\phi$ is the identity outside of $A_0^{m_0}$ and $\eta$ fixes the set of $m_0$ blocks and $m_0$ locations.) 
\end{enumerate}
By items (\ref{I:1}) and (\ref{I:2}) $\phi^{-1} \eta \phi \in \mathcal{G}_{m_0}$.  Since $m_0 > j_0$, $\phi^{-1} \eta \phi  \in \mathcal{G}_{j_0} \subseteq \Sym_{\mathcal{G}}(v_0)$.  It follows that $\eta (\phi(v_0)) = \phi(v_0)$ so $\eta(x) = x$.
\par 
Secondly it is straightforward to prove that the set of pairs  $\{ (\phi^* \mathcal{K},\phi(v_0)) : \phi \in \mathcal{G}(n_0, m_0) \}$ is a one to one function from $\mathcal{H} \slash \mathcal{K}$ onto $V'$.  Since $| \mathcal{H} \slash \mathcal{K}| > \aleph_0$, $|V'| > \aleph_0$.  Let $E' = \{ \{ v_1, v_2 \} : \{ v_1, v_2 \} \subseteq V' \land \{ v_1, v_2 \} \in E \}$.  Since $V'$ is well orderable in $\mathcal{N}9$, using Lemma \ref{L:BG'inFMmodels}, we conclude that either $V'$ has an infinite subset $W$ such that $W$ is independent in $G' = (V', E')$ or $G'$ has a subgraph $G'' = (V'', E'')$ such that $|V''| > \aleph_0$ and $G''$ is a clique in $G'$.  But an independent subset of $V'$ in $G'$ is independent in $G$ and a clique in $G'$ is a clique in $G$.  Therefore, either $V$ has an infinite subset which is independent in $G$ or a subgraph $G'' = (V'', E'')$ such that $G''$ is a clique in $G$ 

\par
\vskip.05in\noindent\textbf{Case 2.} Since $\mathcal{H} \slash \mathcal{K}$ is well-orderable in $\mathcal{N}9$ the only other possible case is $|\mathcal{H} \slash \mathcal{K}| \le \aleph_0$.  In this case we apply Theorem \ref{T:DNT} with $X = T$ and $\mathcal{L} = \mathcal{K}$.  So there is a finite subset $S_1$ of $T$ such that 
\begin{equation} \label{E:S_1_2property}
\{ \eta \in \Sym(T): \forall B \in S_1, \eta(B) = B \} \le \mathcal{K} \le \{ \eta \in \Sym(T) : \eta[S_1] = S_1 \}. 
\end{equation}
\par 
\begin{lemma}
\label{L:S_1ne0}
$S_1 \ne \emptyset$.
\end{lemma}
\begin{proof}
Since $\alpha \in \mathcal{G}_{n_0}$ and $\sup(\alpha) \subseteq A_0^{i_0} \subseteq A_0^{m_0}$ we have $\alpha \in \mathcal{G}(n_0, m_0)$.  Therefore $\alpha^* \in \mathcal{H} = \Sym(T)$.  Since $\alpha(v_0) \ne v_0$, $\alpha^* \notin \mathcal{K}$ and therefore, by the first $\le$ in equation (\ref{E:S_1_2property}, $\alpha^* \notin \{ \eta \in \Sym(T) : \forall s \in S_1, \eta(s) = s \}$.  So there is an $B \in S_1$ such that $\alpha^*(B) \ne B$.
\end{proof}

\begin{lemma}
\label{L:S_1PropGeneralized}
\[
\{ \eta \in \mathcal{G}_{n_0} : \forall B \in S_1, \eta(B) = B \} \le \{ \phi \in \mathcal{G}_{n_0} : \phi(v_0) = v_0 \} \le \{ \eta \in \mathcal{G}_{n_0} : \forall B \in S_1, \eta(B) \in S_1 \}.
\]
\end{lemma}
\begin{proof}
Assume $\gamma \in \mathcal{G}_{n_0}$.  

\begin{sublemma}
\label{SL:betaProps}
$\exists \beta \in \mathcal{G}_{j_0}$ such that 
\begin{enumerate}[(a)]
\item \label{CI:bebInG} $\beta \gamma \beta^{-1} \in \mathcal{G}(n_0, m_0)$ and
\item \label{CI:bebB=B} $\forall B \in S_1$, $\beta(B) = \beta^{-1}(B) = B$.
\item \label{CI:betav_0=v_0} $\beta(v_0) = \beta^{-1}(v_0) = v_0$.
\end{enumerate}
\end{sublemma}
\begin{proof}
Since $\gamma \in \mathcal{G}$ there is a positive integer $k > m_0$ such that $\sup(\gamma) \subseteq A_0^{k}$.  Let $k_0 = k + 1$.  Let
\begin{align*}
W_1 &= \{ C \in \mathcal{B}(j_0, k_0) : C \subseteq A_0^{k_0} \setm A_0^k \} \\
W_2 &= \{ C \in \mathcal{B}(j_0, k_0) : C \subseteq A_0^k \setm A_0^{m_0} \} \mbox{ and } \\
W_3 &= \left\{ C \in \mathcal{B}(j_0, k_0) : C \subseteq A_0^{m_0} \setm (A_0^{j_0} \cup (\bigcup S_1)) \right\}.
\end{align*}
The sets $W_1$, $W_2$ and $W_3$ are all countably infinite.  They are also pairwise disjoint.  It follows that there is a one to one function $F$ from $W_1 \cup W_2 \cup W_3$ onto $W_1 \cup W_2 \cup W_3$ such that 
\begin{equation}
\label{E:propsOfF}
F[W_1] = W_1 \cup W_2, \hskip.1in F[W_2] \subseteq W_3 \hskip.1in \mbox{ and } \hskip.1in F[W_3] \subseteq W_3.
\end{equation}
Let $\beta$ be the element of $\mathcal{G}_{j_0}$ for which 
\[
\forall C \in W_1 \cup W_2 \cup W_3, \beta[C] = F(C) \mbox{ and } \forall C \notin W_1 \cup W_2 \cup W_3, \beta[C] = C.
\]
\par
Since $\beta \in \mathcal{G}_{j_0}$, $\beta$ and $\beta^{-1}$ are in $\mathcal{G}_{n_0}$.  Therefore, since $\gamma$ is also in $\mathcal{G}_{n_0}$ we have $\beta \gamma \beta^{-1} \in \mathcal{G}_{n_0}$.  So to complete the proof of \ref{CI:bebInG} we have to show that $\forall s \in A \setm A_0^{m_0}, \beta \gamma \beta^{-1}(s) = s$.  Assuming $s \notin A_0^{m_0}$ we have $\beta^{-1}(s) \in \bigcup W_1$.  So, since $\bigcup W_1 \cap \sup(\gamma) = \emptyset$,  $\gamma(\beta^{-1}(s)) = \beta^{-1}(s)$ and it follows that $\beta \gamma \beta^{-1}(s) = s$.  
\par 
For part \ref{CI:bebB=B} assume $B \in S_1$.  Every $j_0$ block in $W_1 \cup W_2 \cup W_3$ is disjoint from $\bigcup S_1$ so the $j_0$ block containing $B$ is not in $W_1 \cup W_2 \cup W_3$.  Therefore $\beta(B) = B$ and $\beta^{-1}(B) = B$.  
\par
Part \ref{CI:betav_0=v_0} of the Sublemma holds because $\beta \in \mathcal{G}_{j_0} \subseteq \Sym_{G}(v_0)$. This completes the proof of the sublemma.
\end{proof}

\noindent \textbf{The proof of the first $\le$ in the lemma:}
\par
In addition to the assumption that $\gamma \in \mathcal{G}_{n_0}$, assume that $\forall B \in S_1$, $\gamma(B) = B$.  We have to show that $\gamma(v_0) = v_0$.  
\par 
By our assumption $\gamma(B) = B$ and part \ref{CI:bebB=B} of Sublemma \ref{SL:betaProps},  we conclude that $\beta \gamma \beta^{-1}(B) = B$.
Combining this with part \ref{CI:bebInG} of the sublemma we get
\[
(\beta \gamma \beta^{-1})^* \in \{ \eta \in\Sym(T) : \forall B \in S_1, \eta(B) = B \}.
\]
So by equation (\ref{E:S_1_2property}), $(\beta \gamma \beta^{-1})^* \in \mathcal{K}$ and hence $\beta \gamma \beta^{-1} (v_0) = v_0$.  Using part \ref{CI:betav_0=v_0} of the sublemma the previous equation gives us $\gamma(v_0) = v_0$. 

\noindent \textbf{The proof of the second $\le$ in the lemma:}
\par
In addition to the assumption that $\gamma \in \mathcal{G}_{n_0}$, assume that $\gamma(v_0) = v_0$.  We want to show that $\forall B \in S_1, \gamma(B) \in S_1$.  By Sublemma \ref{SL:betaProps}, part \ref{CI:bebInG}, $\beta \gamma \beta^{-1} \in \mathcal{G}(n_0, m_0)$ and by part \ref{CI:betav_0=v_0} and our assumption, $\beta \gamma \beta^{-1}(v_0) = v_0$.  Therefore $(\beta \gamma \beta^{-1})^* \in \mathcal{K}$ and by (\ref{E:S_1_2property}), $(\beta \gamma \beta^{-1})^*[S_1] = S_1$.  It follows that $\beta \gamma \beta^{-1}[S_1] = S_1$, completing the proof of the second $\le$.
\end{proof} 
\par 
\begin{corollary}
\label{C:phiPsiAgree}
If $\phi$ and $\psi$ are in $\mathcal{G}_{n_0}$ and $\forall C \in S_1, \phi(C) = \psi(C)$ then $\phi(v_0) = \psi(v_0)$.
\end{corollary}

\begin{definition}
\label{D:beta_CD} Assume that $C$ and $D$ are in $\mathcal{B}^{n_0} \setminus \{ A_0^{n_0} \}$ and $C \ne D$ then $\beta_{(C,D)}$ is the element of $\mathcal{G}_{n_0}$ for which $\beta_{(C,D)} \restriction \mathcal{B}^{n_0}$ is the transposition $(C,D)$.  That is, if $s \in A$ then
\[
\beta_{(C,D)}(s) =\begin{cases}
                  (s_0, s_1, \ldots, s_{n_0 - 1})^{\frown} bc^{n_0}(D) & \mbox{ if } s \in C \\
                  (s_0, s_1, \ldots, s_{n_0 - 1})^{\frown} bc^{n_0}(C) & \mbox{ if } s \in D \\
                  s & \mbox{ otherwise}
                  \end{cases} 
\]
\end{definition}
Choose an element $C_0 \in S_1$ (Lemma \ref{L:S_1ne0}) and let 
\begin{equation}
V' = 
\left\{ \beta_{(C_0,C)}(v_0) : C \in \mathcal{B}^{n_0} \setm (S_1 \cup \{ A_0^{n_0} \}) \right\}.
\end{equation}
$V' \subseteq V$ since for every $C \in \mathcal{B}^{n_0} \setm (S_1 \cup \{ A_0^{n_0} \} )$, $\beta_{(C_0, C)} \in \mathcal{G}_{n_0} \subseteq \Sym_{\mathcal{G}} (G)$.
\par 
Define the function $H: \mathcal{B}^{n_0} \setm (S_1 \cup \{ A_0^{n_0} \}) \to V'$ by
\[
H(C) = \beta_{(C_0,C)}(v_0).
\]
\begin{lemma} 
\label{L:propsOfH}
\begin{enumerate}
\item \label{LI:HinN9} $H \in \mathcal{N}9$.
\item \label{LI:H1-1} $H$ is one to one.
\item \label{LI:HontoV'} The range of $H$ is $V'$.
\end{enumerate}
\end{lemma}
\begin{proof}
To prove item (\ref{LI:HinN9}) we show that $\mathcal{G}_{m_0} \subseteq \Sym_{\mathcal{G}}(H)$.  Assume $\phi \in \mathcal{G}_{m_0}$ then $\phi \in \mathcal{G}_{n_0}$ and $\forall B \in S_1$, $\phi(B) = B$.  So $\phi(\dom(H)) = \dom(H)$.  In addition, if $E \in \dom(H)$ then for all $C \in S_1$,
\[
\phi \beta_{(C_0,E)} (C) = \begin{cases}
                             B & \mbox{ if } B \ne C_0 \\
                             \phi(E) & \mbox{ if } B = C_0.
                           \end{cases}
\]
So for all $C \in S_1$, $\phi \beta_{(C_0, E)}(C) = \beta_{(C_0, \phi(E))} (C)$ which, by Corollary \ref{C:phiPsiAgree}, implies that $\phi \beta_{(C_0, E)}(v_0) = \beta_{(C_0, \phi(E))}(v_0)$.  Using the definition of $H$ this equation is equivalent to $\phi(H(E)) = H(\phi(E))$.  Therefore by Lemma \ref{L:WOSetsInFMModels} Part (\ref{L:FMModsPt3}) $\phi(H) = H$.
\par 
For part (\ref{LI:H1-1}) of the Lemma assume that $H(C) = H(D)$.  Then 
\[
\beta_{(C_0, C)}(v_0) = \beta_{(C_0, D)}(v_0) \mbox{ so }\beta_{(C_0, D)}^{-1} \beta_{(C_), C)}(v_0) = v_0.
\]  
Hence by the second $\le$ in Lemma \ref{L:S_1PropGeneralized}, $\forall B \in S_1, \beta_{(C_0, D)}^{-1} \beta_{(C_), C)}(B) = B$.  Letting $B = C_0$ we get $\beta_{(C_0, C)}(C_0) = \beta_{(C_0, D)}(C_0)$ which is equivalent to $C = D$.
\par 
Part \ref{LI:HontoV'} is clear.
\end{proof}

\begin{corollary}
\label{C:PropsOfV'}
$V' \in \mathcal{N}9$ and $|V'| > \aleph_0$ in $\mathcal{N}9$.
\end{corollary}
\begin{proof}
Since $H$ is a one to one correspondence from $\mathcal{B}^{n_0} \setm (S_1 \cup \{ A_0^{n_0} \})$ to $V'$ and $H \in \mathcal{N}9$, $|V'| = |\mathcal{B}^{n_0} \setm (S_1 \cup \{A_0^{n_0}\} )| > \aleph_0$ in $\mathcal{N}9$.
\end{proof}

\begin{lemma}
\label{L:matchingPairs}
Assume $\{ u_1, u_2, v_1, v_2\} \subseteq V'$, $u_1 \ne u_2$ and $v_1 \ne v_2$.  Then there is an $\eta \in \mathcal{G}_{n_0}$ such that $\eta(u_1) = v_1$ and $\eta(u_2) = v_2$.
\end{lemma}
\begin{proof}
We first prove the lemma assuming that $\{ u_1, u_2 \} \cap \{ v_1, v_2 \} = \emptyset$.  Assume $u_1 = \beta_{(C_0, C_1)}(v_0), u_2 = \beta_{(C_0, C_2)}(v_0), v_1 = \beta_{(C_0, D_1)}(v_0)$ and $v_2 = \beta_{(C_0, D_2)}(v_0)$ where $\{ C_1, C_2, D_1, D_2 \} \subseteq \mathcal{B}^{n_0} \setm (S_1 \cup \{ A_0^{n_0} \})$.  Let $\eta = \beta_{(C_1, D_1)} \beta_{(C_2, D_2)}$.  Then $\eta \beta_{(C_0, C_1)} (C_0) = D_1 = \beta_{(C_0, D_1)}(D_1)$ and, since $\{C_1, C_2, D_1, D_2 \} \cap S_1 = \emptyset$, for every $C \in S_1 \setm \{ C_0 \}$, $\eta \beta_{(C_0, C_1)} (C) = C = \beta_{(C_0, D_1)}(C)$. Therefore, by Corollary \ref{C:phiPsiAgree}, $\eta \beta_{(C_0, C_1)} (v_0) = \beta_{(C_0, D_1)} (v_0)$.  So $\eta(u_1) = v_1$.  Similarly, $\eta(u_2) = v_2$.
\par 
For the proof in the general case (dropping the assumption that $\{ v_1, v_2 \} \cap \{u_1, u_2 \} = \emptyset$) we choose $w_1$ and $w_2$ in $V'$ so that $\{ w_1, w_2 \} \cap \{ v_1, v_2, u_1, u_2 \} = \emptyset$.  By the result of the paragraph above there are $\tau$ and $\gamma$ in $\mathcal{G}_{n_0}$ for which $\tau(u_1) = w_1, \tau(u_2) = w_2, \gamma(w_1) = v_1$ and $\gamma(w_2) = v_2$. Then $\eta = \gamma \tau$ satisfies the conclusion of the lemma. 
\end{proof}
\par
We can now complete the proof of Theorem \ref{T:BGinN9} in Case 2.  Choose a pair $\{ v_1, v_2 \} \subseteq V'$ with $v_1 \ne v_2$.  There are two possibilities:  Either $\{v_1, v_2 \} \in E$ or $\{v_1, v_2 \} \notin E$.  In the first case for any two element subset $\{ u_1, u_2 \}$ of $V'$, by Lemma \ref{L:matchingPairs}, there is an $\eta \in \mathcal{G}_{n_0}$ such that $\eta(\{ v_1, v_2 \}) = \{ u_1, u_2 \}$.  Since $\eta \in \Sym_{\mathcal{G}}(G)$, $\{ u_1, u_2 \} \in E$.  So the graph $G' = (V', E')$ where $E' = \{ \{ u_1, u_2 \} \in E : \{u_1, u_2 \} \subseteq V' \}$ is a clique in $G$ with $|V'| > \aleph_0$.  Similarly in the second case $V'$ is an independent subset of $V$.  In either case the conclusion of $\BG$ holds.
\end{proof}

\subsection{The Model $\mathcal{N}26$}\hfill\\ \vskip-.1in

This is the Brunner-Pincus model which is described in \cite{Howard-Rubin-1998}. The set of atoms $A=\bigcup_{n\in\omega} P_n$, where the $P_n$'s are pairwise
disjoint denumerable sets; $\mathcal{G}$  is the group of all permutations
$\sigma$ of $A$ such that $\sigma(P_n)=P_n$, for all $n\in\omega$; and the normal filter $\Gamma = \{ H \le \mathcal{G} : \mbox{ for some finite } E \subseteq A, \fix_{\mathcal{G}}(E) \subseteq H \}$.
A.\ Banerjee has pointed out that 
\begin{enumerate}
\item $\T$ is true in the model as proved by E.\ Tachtsis in \cite[Remark 3, part 2]{Tachtsis-2024a}. 
\item $\BG$ is false in the model since 
  \begin{enumerate}
  \item $\BG \Rightarrow \mathsf{DF{=}F}$ by Proposition \ref{P:impliesDF=F} and the Definition of $\BG$, 
  \item $\mathsf{DF{=}F}$ is known to be false in the model (see \cite{Howard-Rubin-1998}). 
  \end{enumerate}
\end{enumerate}
This gives a complete answer to Question 6.5 in \cite{Banerjee-Gopaulsingh-2023}.

\subsection{The Model $\mathcal{N}17$}\label{SS:Model_N17}\hfill\\ \vskip-.1in

Recall that, by Lemma \ref{lem:EDMT_CCC}, $\T$ (and thus $\BG$) implies $\mathsf{AC}_{\aleph_{0}}^{\aleph_{0}}$. This yields the following natural question: Is there an uncountable well-ordered cardinal $\kappa$ such that $\BG$ (or $\T$) implies $\mathsf{AC}_{\kappa}^{\aleph_{0}}$? This is equivalent to asking whether or not $\BG$ (or $\T$) implies $\mathsf{AC}_{\aleph_{1}}^{\aleph_{0}}$. The purpose of this subsection is to provide a negative answer to the latter question by using the model $\mathcal{N}17$ from \cite{Howard-Rubin-1998}. This will also give us that, for any uncountable well-ordered cardinal $\kappa$, $\BG$ (and thus $\T$) does not imply $\mathsf{UT}(\aleph_{0},\kappa,\kappa)$ in $\mathsf{ZFA}$. However, whether or not $\BG$ (or $\T$) implies $\mathsf{UT}(\aleph_{0},\aleph_{0},\aleph_{0})$ is still an \emph{open problem}.
\par
Model $\mathcal{N}17$ originally appeared in Brunner and Howard \cite{Brunner-Howard-1992} as a member of a class of permutation models defined therein. The description of $\mathcal{N}17$ that is given below, is, on the one hand, different from, but equivalent to, the description given in \cite{Brunner-Howard-1992} and, on the other hand, a correction of the corresponding one in \cite{Howard-Rubin-1998}. 
\par
We start with a model $M$ of $\mathsf{ZFA}+\mathsf{AC}$ with a set $A$ of atoms which is a denumerable disjoint union $A=\bigcup_{n\in\omega}A_{n}$, where $|A_{n}|=\aleph_{1}$ for all $n\in\omega$. Let 
$$G=\{\phi\in\Sym(A):(\forall n\in\omega)(\phi(A_{n})=A_{n})\wedge(\exists J\subseteq\omega)(|J|<\aleph_{0}\text{ and $\phi$ is the identity map on $\bigcup_{i\in\omega\setminus J}A_{i}$})\}.$$ 
Let $\mathcal{F}$ be the filter of subgroups of $G$ generated by the subgroups $\fix_{G}(E)$, where $E=\bigcup_{i\in I}A_{i}$ for some finite $I\subseteq \omega$; $\mathcal{F}$ is a normal filter on $G$. 
$\mathcal{N}17$ is the permutation model determined by $M$, $G$ and $\mathcal{F}$. 
\par
For every $x\in\mathcal{N}17$, there is a finite $I\subseteq\omega$ such that $\fix_{G}(\bigcup_{i\in I}A_{i})\subseteq\Sym_{G}(x)$; any such finite union $\bigcup_{i\in I}A_{i}$ will be called a \emph{support} of $x$.
\par
The subsequent group-theoretic result due to Gaughan \cite{Gaughan-1964} will be useful for the proof of Theorem \ref{thm:BG_UT} below. (Also see Remark \ref{rem:conclude}(\ref{Q4}) below.)

\begin{theorem}
\label{thm:G}
[\cite{Gaughan-1964}] Let $\kappa,\lambda$ be two well-ordered cardinals such that $\aleph_{0}<\kappa\leq\lambda^{+}$. Then the group
$$S(\lambda,\kappa)=\{\sigma\in\Sym(\lambda):|\sup(\sigma)|<\kappa\}$$
has no proper subgroups of index less than $\lambda$. 
\end{theorem}  

\begin{theorem}
\label{thm:BG_UT}
$\BG$ (and thus $\T$) does not imply $\mathsf{AC}_{\aleph_{1}}^{\aleph_{0}}$ in $\mathsf{ZFA}$. Hence, for any uncountable well-ordered cardinal $\kappa$, $\BG$ does not imply $\mathsf{AC}_{\kappa}^{\aleph_{0}}$ (and thus $\mathsf{UT}(\aleph_{0},\kappa,\kappa)$) in $\mathsf{ZFA}$.
\end{theorem}

\begin{proof}
Considering the model $\mathcal{N}17$, we first observe that, for every $n\in\omega$, $|A_{n}|=\aleph_{1}$ in $\mathcal{N}17$, since $A_{n}$ is a support of every element of $A_{n}$ (and thus $A_{n}$ is well orderable in $\mathcal{N}17$) and $|A_{n}|=\aleph_{1}$ in $M$. Furthermore, the family $\mathcal{A}=\{A_{n}:n\in\omega\}$ is denumerable in $\mathcal{N}17$ (since $\fix_{G}(\mathcal{A})=G\in\mathcal{F}$). As mentioned in \cite{Howard-Rubin-1998},  $\mathsf{AC}_{\aleph_{1}}^{\aleph_{0}}$ fails in $\mathcal{N}17$ for $\mathcal{A}$. (In fact, using standard Fraenkel--Mostowski techniques, it is easy to verify that there are no infinite $\mathcal{B}\subseteq\mathcal{A}$ and a function $f$ in $\mathcal{N}17$ such that, for every $B\in\mathcal{B}$, $f(B)$ is a non-empty proper subset of $B$.)

To show that $\mathsf{EDM_{BG}}$ is true in $\mathcal{N}17$, we first prove the following lemma.

\begin{lemma}
\label{lem:1}
In $\mathcal{N}17$, every non-well-orderable set has an uncountable well-orderable subset.
\end{lemma}

\begin{proof}
Let $x\in\mathcal{N}17$ be a set which is not well orderable in $\mathcal{N}17$. Since $x\in\mathcal{N}17$, there exists a finite $K\subseteq\omega$ such that $E=\bigcup_{k\in K}A_{k}$ is a support of $x$. As $x$ is not well orderable in $\mathcal{N}17$, there exists $y\in x$ such that $E$ is not a support of $y$. Let $K_{0}$ be a finite subset of $\omega$ which is disjoint from $K$ and is such that the set $E\cup E_{0}$, where $E_{0}:=\bigcup_{k\in K_{0}}A_{k}$, is a support of $y$. 

Since $E$ is not a support of $y$, there exists $\phi\in\fix_{G}(E)\setm\Sym_{G}(y)$. Let $\eta$ be the permutation of $A$ which agrees with $\phi$ on $E_{0}$ and $\eta$ is the identity map on $A\setminus E_{0}$. Clearly, $\eta\in G$. Moreover, $\eta\in\fix_{G}(E)$, since $E\cap E_{0}=\emptyset$. As $\eta$ agrees with $\phi$ on $E\cup E_{0}$ and $E\cup E_{0}$ is a support of $y$, we have $\eta(y)=\phi(y)$, and thus $\eta(y)\neq y$, since $\phi(y)\neq y$. Furthermore, as $K_{0}$ is finite and $\eta=\prod_{k\in K_{0}}\eta_{k}$, where for $k\in K_{0}$, $\eta_{k}\upharpoonright A_{k}=\phi\upharpoonright A_{k}$ and $\eta_{k}\upharpoonright (A\setm A_{k})$ is the identity map, it follows that for some $k_{0}\in K_{0}$, $\eta_{k_{0}}(y)\neq y$.
We let $$\mathbf{G}=\fix_{G}(A\setm A_{k_{0}}).$$ 
Then $\mathbf{G}$ is isomorphic to the group $S(\aleph_{1},\aleph_{2})$, which is isomorphic to $\Sym(\aleph_{1})$. Furthermore,  $\Orb_{\mathbf{G}}(y)$ ($=\{\pi(y):\pi\in \mathbf{G}\}$) is a subset of $x$ (since $y\in x$ and $\mathbf{G}\subseteq\fix_{G}(E)\subseteq\Sym_{G}(x)$) which is well orderable in $\mathcal{N}17$. Indeed, $\fix_{G}(E\cup E_{0})\subseteq\fix_{G}(\Orb_{\mathbf{G}}(y))$. To see this, let $\rho\in\fix_{G}(E\cup E_{0})$ and $\pi\in\mathbf{G}$. By definition of $G$, it follows that $\pi(E\cup E_{0})=E\cup E_{0}$. Furthermore, since $E\cup E_{0}$ is a support of $y$, $\pi(E\cup E_{0})$ is a support of $\pi(y)$ and, since $\pi(E\cup E_{0})=E\cup E_{0}$, we get that $E\cup E_{0}$ is a support of $\pi(y)$. Therefore, $\rho(\pi(y))=\pi(y)$, i.e. $\rho\in\fix_{G}(\Orb_{\mathbf{G}}(y))$, as required.\footnote{The latter argument implicitly shows that more is true, namely, for any $z\in\mathcal{N}17$, $\Orb_{G}(z)$ is well orderable in $\mathcal{N}17$.}
\par
To complete the proof of the lemma, we will show that $\Orb_{\mathbf{G}}(y)$ is uncountable in $\mathcal{N}17$. Let
$$\mathbf{H}=\{\pi\in\mathbf{G}:\pi(y)=y\}.$$
Then $\mathbf{H}$ is a proper subgroup of $\mathbf{G}$ (since $\eta_{k_{0}}\in\mathbf{G}\setminus\mathbf{H}$) and $|\Orb_{\mathbf{G}}(y)|=(\mathbf{G}:\mathbf{H})$. As $\mathbf{G}$ is isomorphic to $S(\aleph_{1},\aleph_{2})$ and $\mathbf{H}$ is a proper subgroup of $\mathbf{G}$, it follows from Theorem \ref{thm:G} that $(\mathbf{G}:\mathbf{H})\not<\aleph_{1}$. Therefore $|\Orb_{\mathbf{G}}(y)|\not\leq\aleph_{0}$, as required.         
\end{proof}
\par
From Lemmas \ref{L:BG'inFMmodels} and \ref{lem:1}, we now readily obtain that $\BG$ is true in $\mathcal{N}17$, finishing the proof of the theorem.
\end{proof}

\begin{remark}
\label{UT:Model_N17}
\begin{enumerate}
\item \label{rem:UT_N17} $\mathsf{UT}(\aleph_{0},\aleph_{0},\aleph_{0})$ is true in the model $\mathcal{N}17$, as shown in \cite{Brunner-Howard-1992}. In view of the proof of Lemma \ref{lem:1} and for the reader's convenience, we include the short argument. Let $\mathcal{U}=\{U_{n}:n\in\omega\}$ be a denumerable family of denumerable sets in $\mathcal{N}17$; the mapping $n\mapsto U_{n}$ is a bijection in $\mathcal{N}17$. Let $E$ be a support of the mapping $n\mapsto U_{n}$, $n\in\omega$. Then $E$ is a support of $U_{n}$ for all $n\in\omega$.

We assert that $E$ is a support of every element of $\bigcup\mathcal{U}$. This will give us that $\bigcup\mathcal{U}$ is well orderable in $\mathcal{N}17$, and therefore, since $\bigcup\mathcal{U}$ is denumerable in the ground model $M$, $\bigcup\mathcal{U}$ is denumerable in $\mathcal{N}17$.

Aiming for a contradiction, suppose that, for some $n_{0}\in\omega$, there exists $y\in U_{n_0}$ such that $E$ is not a support of $y$. Now, working almost identically to the proof of Lemma \ref{lem:1}, it can be shown that the orbit of $y$ with respect to some suitable subgroup of $G$ is an uncountable (well-orderable) subset of $U_{n_0}$ in $\mathcal{N}17$. But this contradicts the fact that $U_{n_{0}}$ is denumerable in $\mathcal{N}17$. Therefore $\mathsf{UT}(\aleph_{0},\aleph_{0},\aleph_{0})$ is true in  $\mathcal{N}17$ as required.

\item \label{rem:ET} Tachtsis \cite{Tachtsis-2016} considered the following variant of $\mathcal{N}17$. The set $A$ of atoms and the normal filter $\mathcal{F}$ are defined as in the description of $\mathcal{N}17$, but the group of permutations of $A$ is given by $$G^{*}:=\{\phi\in\Sym(A):(\forall n\in\omega)(\phi(A_{n})=A_{n})\}.$$

Denoting the resulting model by $\mathcal{N}$, it was shown in \cite{Tachtsis-2016} that $\mathsf{DF=F}\wedge\mathsf{UT}(\aleph_0,\aleph_0,\aleph_0)$ is true in $\mathcal{N}$. However, we would like to point out that $\mathcal{N}=\mathcal{N}17$. A proof of this fact can be given as in Howard and Tachtsis \cite[Proof of Theorem 11]{Howard-Tachtsis-2023} by making the obvious modifications therein. More specifically, replacing ``$A_{n}$ is denumerable'' in the statement of \cite[Theorem 11]{Howard-Tachtsis-2023} by ``$|A_{n}|=\aleph_{1}$'', the proof of \cite[Theorem 11]{Howard-Tachtsis-2023} works almost unchanged to show that $\mathcal{N}=\mathcal{N}17$. We thus refer the reader to \cite{Howard-Tachtsis-2023} for the details. Lastly, let us note that, since $\mathsf{DF=F}$ $\wedge$ ``There exists a Dedekind infinite set without uncountable well-orderable subsets'' is relatively consistent with $\mathsf{ZF}$ (by using, for example, $\mathcal{N}_{T}$ and a transfer theorem of Pincus), it follows that Lemma \ref{lem:1} is a proper strengthening of `$\mathcal{N}\models\mathsf{DF=F}$' in \cite{Tachtsis-2016}.      
\end{enumerate}   
\end{remark}

\section{Summary of the results} \label{S:Diagram}
Our results are summarized by the following diagram.  The subscript $\mathsf{tr}$ on a Fraenkel-Mostowski model name means that the result transfers to $\mathsf{ZF}$.
\vskip.2in

\begin{center}
\begin{tikzpicture}
\node (DM) at (-4,4) {$\boxed{\DM \leftrightarrow \mathsf{AC}}$};
\node (BG) at (-4,2) {$\BG$};

\node (AC^om_om1) at (-6.5,0.5) {$\mathsf{AC}_{\aleph_{1}}^{\aleph_{0}}$};

\node (2m=m) at (-.5,3) {$2m{=}m$};

\node (PC2and4) at (4,2.5) {$\boxed{\mathsf{PC}(\aleph_0, 2, \aleph_0)  \wedge \mathsf{PC}(\aleph_0, 4, \aleph_0)} $};

\node (T) at (-4,0) {$\T$};
\node (DF=F) at (-.5,0) {$\mathsf{DF{=}F}$};
\node (PC>) at (2.5,-2.5) {$\boxed{\mathsf{PC}({>} \aleph_0, {<}\aleph_0, {>}\aleph_0) \leftrightarrow (\mathsf{PC}({>} \aleph_0, {<}\aleph_0, {\nleq}\aleph_0) \wedge \mathsf{AC_{fin}^{\aleph_0}})} $};

\node(PCnleq2) at (2.5,-6) {$\mathsf{PC}({\nleq}\aleph_0,{<}\aleph_{0},{>}\aleph_0)$};

\node (ACa_0a_0) at (-3.5,-2.5) {$\mathsf{AC}_{\aleph_0}^{\aleph_0}$};
\node (PCnleq) at (6,0) 
{$\boxed{\mathsf{PC}({>} \aleph_0, {<}\aleph_0, {\nleq}\aleph_0)\leftrightarrow\mathsf{PC}({\nleq} \aleph_0, {<}\aleph_0, {\nleq}\aleph_0)}$};

\draw [thick, ->] (ACa_0a_0) -- (PC>) node[pos=0.3,sloped] { $/$} node[pos=0.5, below] {$\mathcal{N}_{T}$};

\draw [thick,->] (BG) -- (AC^om_om1) node[pos=0.3,sloped] { $/$} node[pos=0.5, left] {$\mathcal{N}17$};

\draw [thick, <-](PC>.-38) -- (PCnleq2.28);

\draw [thick,->] (PC>) -- (PCnleq2) node[pos=0.3,sloped] { $/$} node[pos=0.5, left] {$\mathcal{N}1, \mathcal{N}3$};

\draw [thick, ->](DM.-40) -- (BG.33);
\draw [thick, ->](BG.150) -- (DM.-140) node[pos=0.3, sloped] {$/$} node[pos=0.5, left] {$\mathcal{N}12(\aleph_1), \mathcal{N}9$};
\draw [thick,->] (BG.20) -- (2m=m) node[pos=0.5,sloped] { $/$} node[pos=0.5, below, sloped] {$\mathcal{N}12(\aleph_1)$};
\draw [thick, ->] (DM) -- (2m=m);
\draw [thick, ->](BG.-90) -- (T.90);
\draw [thick,->] (T.10) -- (DF=F.170) node[pos=0.3,sloped] { $/$} node[pos=0.4, above, sloped] {$\mathcal{N}1, \mathcal{N}3$};
\draw [thick, ->] (BG.-15) -- (DF=F.150);
\draw [thick, ->] (2m=m) -- (DF=F);
\draw [thick,->] (DF=F.5) -- (PCnleq) node[pos=0.3, sloped] { $/$} node[pos=0.5, above, sloped] {$\mathcal{N}_T$};
\draw [thick, ->] (T) -- (ACa_0a_0);
\draw [thick, ->] (T) -- (PC>.160);
\draw [thick, ->] (PC>) -- (PCnleq);
\draw [thick, ->] (T) -- (ACa_0a_0);
\draw [thick,->] (PC>.170) -- (T.-94) node[pos=0.7,sloped] { $/$} node[pos=0.5, below, sloped] {$\mathcal{N}5_{\mathsf{tr}}$};
\draw [thick, ->] (T) -- (ACa_0a_0);
\draw [thick, ->] (PCnleq) -- (PC2and4);
\draw [thick, ->] (DF=F) -- (PC2and4);
\draw [thick, ->] (BG) to [bend right = 88] (PCnleq2);
\draw [thick, ->] (DM) to [bend right = 67] (AC^om_om1);
\end{tikzpicture}
\end{center}

\section{Questions and concluding remarks}
\label{S:Questions}
\begin{enumerate}
\item \label{Q:PC1toPC2} Does $\mathsf{PC}({>}\aleph_0, {<} \aleph_0, {\nleq} \aleph_0)$ imply $\mathsf{PC}({>}\aleph_0, {<} \aleph_0, {>} \aleph_0)$?
\item \label{Q:PC_n-element} For which $n \in \omega$ does $\mathsf{PC}({>}\aleph_0, {<} \aleph_0, {\nleq} \aleph_0)$ imply $\mathsf{PC}(\aleph_0, n, \aleph_0)$?  (The implication holds for $n \in \{2,4\}$.  See Proposition \ref{P:PCnleqToPC2=}.)
\item Does $\BG$, or $\T$, imply $\mathsf{UT}(\aleph_{0},\aleph_{0},\aleph_{0})$?
\item \label{Q:Choice_for_aleph1} Does $\BG$, or $\T$, imply $\mathsf{AC}_{\aleph_{0}}^{\aleph_{1}}$, or $\mathsf{AC}_{n}^{\aleph_{1}}$ for some $n\in\omega\setm\{0,1\}$, or $\mathsf{AC_{fin}^{\aleph_{1}}}$ (Every $\aleph_{1}$-sized set of non-empty finite sets has a choice function)? 
\end{enumerate}

Regarding Question (\ref{Q:Choice_for_aleph1}), we have the following remarks. Especially, (\ref{Q4}) of Remark \ref{rem:conclude} below provides a partial answer to the first part of Question (\ref{Q:Choice_for_aleph1}). 

\begin{remark}
\label{rem:conclude}
\begin{enumerate}
\item $\BG$ does not imply $\mathsf{AC}_{2}^{\aleph_{2}}$ in $\mathsf{ZFA}$, so neither does it imply $\mathsf{AC}_{\lambda}^{\kappa}$ in $\mathsf{ZFA}$, for all well-ordered cardinals $\kappa\ge\aleph_{2}$ and $\lambda\ge\aleph_{0}$. Indeed, Jech in \cite[Proof of Theorem 8.3]{Jech} has constructed a permutation model $\mathscr{V}$ in which ``For every set $X$, $|X|\leq\aleph_{1}$ or $|X|\geq \aleph_{1}$'' is true, but $\mathsf{AC}_{2}^{\aleph_{2}}$ is false. By the fact that the former statement is true in $\mathscr{V}$ and by Lemma \ref{L:BG'inFMmodels}, it follows that $\BG$ is true in $\mathscr{V}$. Furthermore, a straightforward modification in the construction of Jech's model $\mathscr{V}$ yields $\BG$ does not imply $\mathsf{AC}_{n}^{\kappa}$ in $\mathsf{ZFA}$, for all well-ordered cardinals $\kappa\ge\aleph_{2}$ and for all $n\in\omega\setm\{0,1,2\}$.

\item Note that, by Proposition \ref{P:ConsequencesOfEDM}(\ref{PP:BGtoPC}), $\T$ implies the weak choice principle ``Every $\aleph_{1}$-sized set of non-empty finite sets has an $\aleph_{1}$-sized subset with a choice function''.  

\item \label{Q4} $\BG$ does not imply $\mathsf{AC}_{\ge\aleph_{0}}^{\aleph_{1}}$ (Every $\aleph_{1}$-sized set of Dedekind infinite sets has a choice function) in $\mathsf{ZFA}$. To see this, first observe that
\begin{equation}
\label{eq:rem3}
\mathsf{AC}_{\ge\aleph_{0}}^{\aleph_{1}}\Rightarrow\mathsf{AC}_{\aleph_{1}}^{\aleph_{1}}\Rightarrow\mathsf{AC}_{\aleph_{1}}^{\aleph_{0}}.
\end{equation}
By Theorem \ref{thm:BG_UT}, we know that $\BG$ does not imply $\mathsf{AC}_{\aleph_{1}}^{\aleph_{0}}$ in $\mathsf{ZFA}$. Therefore, by (\ref{eq:rem3}), $\BG$ does not imply $\mathsf{AC}_{\aleph_{1}}^{\aleph_{1}}$, and thus $\mathsf{AC}_{\ge\aleph_{0}}^{\aleph_{1}}$, in $\mathsf{ZFA}$ either.
\par
Finally, let us point out that none of the implications in (\ref{eq:rem3}) are reversible in $\mathsf{ZFA}$. 
\par
That the first implication in (\ref{eq:rem3}) is not reversible in $\mathsf{ZFA}$ can be established by using a variant of Model $\mathcal{N}17$ in Subsection \ref{SS:Model_N17}: The set $A$ of atoms is a disjoint union $A=\bigcup_{\alpha<\omega_{1}}A_{\alpha}$ where, for each $\alpha<\omega_{1}$, $|A_{\alpha}|=\aleph_{2}$. The group $G$ and the normal filter $\mathcal{F}$ on $G$ are defined as in Subsection \ref{SS:Model_N17}, except for replacing $\omega$ by $\omega_{1}$ in their definitions therein. Let $\mathcal{V}$ be the resulting permutation model. Following a line of reasoning similar to the one in the proof of Theorem \ref{thm:BG_UT} and in the argument in Remark \ref{UT:Model_N17}(\ref{rem:UT_N17}), it can be verified that $\mathcal{V}\models\mathsf{UT}(\aleph_{1},\aleph_{1},\aleph_{1})\wedge\neg\mathsf{AC}_{\aleph_{2}}^{\aleph_{1}}$, which yields $\mathsf{AC}_{\aleph_{1}}^{\aleph_{1}}\nRightarrow\mathsf{AC}_{\ge\aleph_{0}}^{\aleph_{1}}$ in $\mathsf{ZFA}$ as required. The details are left as an easy exercise for the reader.
\par
That the second implication in (\ref{eq:rem3}) is not reversible in $\mathsf{ZFA}$ can be verified by using the following permutation model: The set $A$ of atoms is a disjoint union $A=\bigcup_{\alpha<\omega_{1}}A_{\alpha}$ where, for each $\alpha<\omega_{1}$, $|A_{\alpha}|=\aleph_{1}$. The group of permutations of $A$ is $G=\{\phi\in\Sym(A):(\forall\alpha<\omega_{1})(\phi(A_{\alpha})=A_{\alpha})\}$. The ideal of supports is defined by $\mathcal{I}=\{E\in\mathscr{P}(A):\text{ for some countable $F\subseteq\omega_{1}$, }E\subseteq\bigcup_{\alpha\in F}A_{\alpha}\}$. $\mathcal{I}$ is a normal ideal, so $\mathcal{F}=\{H\leq G:(\exists E\in\mathcal{I})(\fix_{G}(E)\subseteq H)\}$ is a normal filter on $G$. Let $\mathscr{U}$ be the permutation model determined by $M$, $G$, and $\mathcal{F}$.

Since $\mathcal{I}$ is closed under countable unions, it follows that the principle of dependent choices (Form 43 in \cite{Howard-Rubin-1998}) is true in $\mathscr{U}$ (see, for example, Jech \cite[Lemma 8.4, proof of Theorem 8.3]{Jech}), and thus $\mathsf{AC}_{\aleph_{1}}^{\aleph_{0}}$ is also true in $\mathscr{U}$. On the other hand, $\mathsf{AC}_{\aleph_{1}}^{\aleph_{1}}$ fails in $\mathscr{U}$ for the collection $\mathcal{A}=\{A_{\alpha}:\alpha<\omega_{1}\}$. Therefore $\mathsf{AC}_{\aleph_{1}}^{\aleph_{0}}\nRightarrow\mathsf{AC}_{\aleph_{1}}^{\aleph_{1}}$ in $\mathsf{ZFA}$ as required.     
\end{enumerate}
\end{remark}

\bibliographystyle{amsalpha}

\section*{Conflict of interest}
 The authors declare that they have no conflict of interest.

\end{document}